\documentclass[11pt,english]{article}
\usepackage{lmodern}
\usepackage[T1]{fontenc}
\usepackage[latin9]{inputenc}
\usepackage{geometry}
\usepackage{color}
\usepackage{babel}
\usepackage{float}
\usepackage{amsmath}
\usepackage{amsthm}
\usepackage{amssymb}
\usepackage{graphicx}
\PassOptionsToPackage{normalem}{ulem}
\usepackage{ulem}
\usepackage[unicode=true]
 {hyperref}

\makeatletter

\providecommand{\tabularnewline}{\\}
\floatstyle{ruled}
\newfloat{algorithm}{tbp}{loa}
\providecommand{\algorithmname}{Algorithm}
\floatname{algorithm}{\protect\algorithmname}

\numberwithin{equation}{section}
\numberwithin{figure}{section}
\theoremstyle{plain}
\newtheorem{thm}{\protect\theoremname}
  \theoremstyle{plain}
  \newtheorem{prop}[thm]{\protect\propositionname}

\usepackage{algorithm}
\usepackage{algorithmic}

\def\erf{\mathop{\operator@font erf}\nolimits}

\makeatother

  \providecommand{\propositionname}{Proposition}
\providecommand{\theoremname}{Theorem}

\begin{document}

\title{Spectral Condition-Number Estimation of Large Sparse Matrices}

\author{{\bf Haim Avron}\\
Tel Aviv University
\and {\bf Alex Druinsky} \\
Tel Aviv University
\and {\bf Sivan Toledo} \\
Tel Aviv University }
\maketitle
\begin{abstract}
We describe a randomized Krylov-subspace method for estimating the
spectral condition number of a real matrix $A$ or indicating that
it is numerically rank deficient. The main difficulty in estimating
the condition number is the estimation of the smallest singular value
$\sigma_{\min}$ of $A$. Our method estimates this value by solving
a consistent linear least-squares problem with a known solution using
a specific Krylov-subspace method called LSQR. In this method, the
forward error tends to concentrate in the direction of a right singular
vector corresponding to $\sigma_{\min}$. Extensive experiments show
that the method is able to estimate well the condition number of a
wide array of matrices. It can sometimes estimate the condition number
when running a dense SVD would be impractical due to the computational
cost or the memory requirements. The method uses very little memory
(it inherits this property from LSQR) and it works equally well on
square and rectangular matrices.
\end{abstract}

\section{Introduction}

This paper discusses the problem of estimating the spectral condition
number of a large and sparse matrix $A\in\mathbb{R}^{m\times n}$.
Without loss of generality we assume that $m\geq n$ (we can estimate
the condition number of $A^{T}$ if $n>m$). The main difficulty is
in estimating the smallest singular value $\sigma_{\min}$ of $A$.
Dense SVD algorithms can approximate $\sigma_{\min}$ well and their
running time is predictable, but they are slow (their running time
depends on the dimensions of the matrix, and not on the number of
non-zeros entries in the matrix). More importantly, dense SVD algorithms
require space that is proportional to $mn$ when the matrix is $m$-by-$n$,
which is impractical for large sparse matrices. Likewise, if $A$
is square and a $LU$ factorization of $A$ can be computed then $\sigma_{\min}$
can be estimated efficiently from the factorization~\cite{GH15},
however for some large sparse matrices space requirements due to fill-in
make this method impractical.

Symmetrization is not always an effective way to address the problem.
If we work with the Gram matrix $A^{T}A$, we cannot estimate condition
numbers $\kappa(A)=\sigma_{\max}/\sigma_{\min}$ greater than $1/\sqrt{\epsilon_{\text{m}}}$,
where $\epsilon_{\text{m}}$ is the unit roundoff (machine precision),
unless the minimum singular value is well separated~\cite{Jia2006}.
If we work with the augmented matrix
\[
\begin{bmatrix}0 & A^{T}\\
A & 0
\end{bmatrix}\;,
\]
$\sigma_{\min}$ is transformed into a pair of eigenvalues $\pm\sigma_{\min}$
in the middle of the spectrum. Such eigenvalues are difficult to compute
accurately with the Lanczos algorithm and its variants\footnote{For example, the ARPACK Users' Guide states that a shift-invert iteration
is usually required to compute eigenvalues in the interior of the
spectrum~\cite[Section 3.4]{ARPACK-UG}. Experiments with ARPACK
on some of the matrices presented later in the paper, whose condition
numbers our method was able to estimate, showed that ARPACK does not
converge on them when it tries to compute the smallest-magnitude eigenvalues
of the augmented matrix without inversion.}, and it is essentially impossible for such algorithms to determine
whether there is an eigenvalue closer to zero than the one that has
already been computed.

If one is interested in the smallest singular value of a large sparse
matrix, then avoiding the impractical space requirements of dense
SVD requires the use of a low-memory iterative method. Some Krylov-subspace
methods estimate the condition number as the iteration progresses
(e.g., LSQR). However, LSQR's condition number estimator has some
shortcomings that make it unreliable. We explain these in the next
section.

The starting point of this paper is the following observation, which
we explain in more detail in Section~\ref{sec:Mathematical-Motivation}.
Suppose $b$ is a vector in the column space of $A$, and let $Ax^{\star}=b$.
Let $x^{(0)},x^{(1)},\dots$ be the iterates of a Krylov-subspace
method applied to the problem $\arg\min\|Ax-b\|$ (all norms in this
paper denote the $2$-norm unless stated otherwise). If we knew the
forward error $d^{(k)}=x^{\star}-x^{(k)}$, we could use the Rayleigh
quotients $\|Ad^{(k)}\|/\|d^{(k)}\|$ as an upper estimate on $\sigma_{\min}$.
Now, if $b$ is also a random vector, then the iterates $x^{(k)}$
tend to generate at some point forward errors $d^{(k)}$ such that
the Rayleigh quotient $\|Ad^{(k)}\|/\|d^{(k)}\|$ estimates $\sigma_{\min}$
well. The underlying reason is the tendency of Krylov-subspace methods
to concentrate the forward error in the direction of a singular vector
associated with $\sigma_{\min}$. This tendency is essentially the
same as the tendency of the Lanczos algorithm to converge to outer
eigenvalues first (see, e.g., ~\cite[Section~7.3]{DemmmelAppliedNumericalLinearAlgebra}).
Both tendencies reflect the uneven convergence of Krylov subspaces
to singular- or eigenspaces, one in the singular unsymmetric case
(LSQR) and the other in the symmetric case (Lanczos or Conjugate Gradients).
This is often seen as a flaw in LSQR and in Conjugate Gradients (LSQR
is mathematically equivalent to Conjugate Gradients applied to $A^{T}Ax=A^{T}b$).
These solvers converge slowly when the coefficient matrix $A$ is
ill conditioned because it is difficult for them to get rid of the
error in the right singular subspaces of $A$ that correspond to its
smallest singular values. However, one can also exploit this flaw,
which acts as a sieve that captures a vector from this subspace.

In this paper, we use this observation to design a Krylov-subspace
method that can estimate $\sigma_{\min}$ (hence the spectral condition
number $\mathbf{\sigma_{\max}/}\sigma_{\min}$) accurately. To have
access to the forward error, our algorithm generates a random $x^{\star}$
and computes $b=Ax^{\star}.$ The algorithm then proceeds with using
LSQR on the problem $\arg\min\|Ax-b\|$, meanwhile monitoring the
value of $\|Ad^{(k)}\|/\|d^{(k)}\|$ and keeping the vector $d^{(\ell)}$
for which it is minimized. Very little memory is needed, as the algorithm's
memory usage is close to LSQR's. As the iteration progresses, the
estimate of $\sigma_{\min}$ does not deteriorate; we always keep
the best $d^{(k)}$ found. However, it is undesirable to continue
to run the algorithm when the estimate of $\sigma_{\min}$ is unlikely
to improve. To that end, we design a set of stopping criteria that
are tuned towards the algorithm's goal of computing the condition
number, as opposed to solving least-squares equations, which is what
LSQR's regular stopping criteria are tuned for. Some of our stopping
criteria are the same as LSQR's, but others are based on the algorithm's
knowledge of the forward error.

Our algorithm will never report an underestimate of $\sigma_{\min}$.
Furthermore, in most cases it also produces a certificate vector that
proves that the estimate of $\sigma_{\min}$ is indeed an upper bound.
While we are unable to prove rigorous bounds on the quality of approximation,
nor a bound on the number of iterations until convergence, extensive
experiments indicate that our method does not incorrectly report significant
overestimates of $\sigma_{\min}$, at least when $A$ is not close
to being numerically rank deficient. When $A$ is close to numerical
rank deficiency, the algorithm reports a condition number near $\epsilon_{\text{m}}^{-1}$,
but the accuracy of the estimate is sometimes poor.

Our method also has some flaws. The main one is that it sometimes
converges very slowly, making it essentially impossible to compute
$\sigma_{\min}$. Our experience shows that the method always converges
eventually, but that convergence might be too slow to be of practical
use. There is no good way to determine how close the method is to
termination, although it tends to behave consistently on related matrices
(e.g., from the same application area). When $\kappa(A)$ is close
to $\epsilon_{\text{m}}^{-1}$, the method sometimes overestimates
$\sigma_{\min}$ by several orders of magnitude. The method sometimes
converges very rapidly and sometimes very slowly. This is not related
to the size of the problem and not to how ill conditioned it is, but
to the distribution of singular values.

Even with these flaws, to the best of our knowledge this method is
the only practical way to compute $\sigma_{\min}$ with reasonable
accuracy (to within a factor of $2$ or better) on many large matrices.
Our implementation of the method in MATLAB is freely available\footnote{See~\href{https://github.com/sparse-condest/condest}{https://github.com/sparse-condest/condest}.},
and a C++ implementation is available inside the libSkylark library\footnote{See~\href{https://xdata-skylark.github.io/libskylark/}{https://xdata-skylark.github.io/libskylark/}.}.

The rest of this paper is organized as follows. Section~\ref{sec:Related-Work}
surveys related work on condition number estimation. Section \ref{sec:Mathematical-Motivation}
presents a mathematical analysis of the quantity that we use to estimate
$\sigma_{\min}$, to motivate the algorithm that we present in Section~\ref{sec:The-Algorithm}.
Section~\ref{sec:analysis-small-err} analyzes a unique stopping
criteria that our method relies upon. Section~\ref{sec:numerical}
illustrates and explores the behavior of our algorithm using numerical
experiments.

\section{\label{sec:Related-Work}Related Work}

The largest singular value $\sigma_{\max}$ of $A$ can be computed
accurately using a bounded number of matrix-vector products involving
$A$ and $A^{T}$. This can be done using the power method, for example,
whose analysis for this application we explain in Section~\ref{sec:The-Algorithm}.
The Lanczos method can reduce the number of matrix-vector products
even further~\cite{Kuczynsky92}. Random projection methods can also
estimate $\sigma_{\max}$~\cite{HalkoMartinssonTropp11}.

Estimating $\sigma_{\min}$ is computationally more challenging, because
applying the pseudoinverse is usually much harder than applying $A$
itself. In general, existing random projection methods cannot efficiently
estimate $\sigma_{\min}$ unless a decomposition of $A$ is computed,
or $A$ is low rank (or numerically low rank). If $A$ is low rank,
random projection methods can be used to estimate the $k$th largest
singular value, where $k$ is the (numerical) rank~\cite{HalkoMartinssonTropp11}.
Another approach is to use iterative methods for computing singular
triplets. Recently, variants of Lanczos that use implicit restarting,
harmonic Ritz values and other advanced techniques to further accelerate
the convergence to the smallest singular triplets have been proposed~\cite{KBG04,BaglamaReichel13}.
Another approach for computing singular triplets (including the smallest
singular triplets) is the JDSVD method \cite{Hochstenbach01} which
extends the Jacobi-Davidson method for singular value problems. The
PRIMME\_SVDS~\cite{PRIMME_SVDS} library implements an hybrid algorithm
PHSVDS~\cite{WS15} which is based on carefully selecting an appropriate
algorithm at each point of the computation.'' Also worth mentioning
is the inverse free preconditioned Krylov subspace method~\cite{GolubYe02,LiangYe14}.

The LINPACK condition-number estimator requires a triangular factorization
of $A$ (see Higham's monograph~\cite[Chapter 15]{Higham02} for
details on this and related estimators). The Gotsman-Toledo~\cite{GostmanToledo2008}
and the Bischof et al.~\cite{Bischof90} condition-number estimators,
which are specialized to sparse matrices, also require a triangular
factorization. Estimators that require a triangular factorization
are less expensive than the SVD, but they still cannot be applied
to huge matrices.

The LAPACK condition-number estimator estimates the 1-norm condition
number~\cite{Higham88} . Higham and Tisseur derive a block version
of the LAPACK 1-norm condition number estimator~\cite{HT00}. Both
methods are black-box in the sense that they require only matrix-vector
products with the pseudoinverses of $A$ and $A^{T}$. One way to
compute these products is using a factorization, but they can also
be computed using an iterative solver. With an effective preconditioner,
repeated applications of the pseudoinverse may be less expensive than
the method that we propose, but without one our method is less expensive.
Kenney et al.~\cite{Kenney98} describe a way to estimate the condition
number of a square matrix using a single application of the inverse
to one or several vectors. Vecharynski proposes a Rayleigh-quotient-type
iteration that uses a single iteration of a preconditioned iterative
solver to compute an approximate singular vector given an approximate
singular value (which is estimated form the previous singular-vector
estimate)~\cite{vecharyn10}.  Gaaf and Hochstenbach~\cite{GH15}
propose to estimate the condition number of square matrices using
the \emph{extended Krylov subspace}. Computing a basis for the extended
Krylov subspace requires applying the inverse of the matrix, so the
method uses a $LU$ factorization of the matrix.

The spectral condition number measures the normwise sensitivity of
matrix-vector products and linear systems to small perturbations in
the inputs. There are methods that estimate more focused metrics,
such as the sensitivity of individual components of the inputs or
output~\cite{Kenney98}. Our method does not address this problem.

LSQR itself also estimates a condition number of $A$, the Frobenius
condition number $\|A\|_{F}\|A^{+}\|_{F}$, where $A^{+}$ is the
pseudoinverse of $A$. LSQR estimates the two norms by computing explicitly
the Frobenius norm of two bidiagonal matrices that it computes as
part of the algorithm, $B^{(t)}$ and $D^{(t)}$. The estimates rely
on inequalities $\|B^{(t)}\|_{F}\leq\|A\|_{F}$ and $\|D^{(t)}\|_{F}\leq\|A^{+}\|_{F}$.
However, this method has several serious shortcomings. First, the
LSQR paper does not provide any evidence that the gap in the inequalities
is small by the time the algorithm terminates. In particular, it appears
difficult to come up with stopping criteria that will ensure that
these estimates are accurate. Second, the inequalities depend on the
theoretical (exact arithmetic) orthogonality of the Lanczos vectors,
which are far from orthogonal when the algorithm is implemented in
floating-point arithmetic.

As noted by Paige and Saunders~\cite{LSQR}, the spectral condition
number can also be bounded using $\|B^{(t)}\|_{2}\|D^{(t)}\|_{2}$
but this is not done in practice (nor is it suggested by Paige and
Saunders) since computing this quantity in each iteration is expensive.

Our method improves upon the condition number estimates of LSQR in
several ways. First, we devise a way to approximate the spectral condition
number cheaply in each iteration, so our algorithm estimates the spectral
condition number instead of the Frobenius condition number, as is
done by LSQR. The spectral condition number is more relevant for estimating
the accuracy of linear solvers. Second, we use two estimators for
$\sigma_{\min}$, one of which is similar in spirit to LSQR's use
of $\|D^{(t)}\|_{F}$ and the other completely novel, which only relies
\emph{on properties that do hold in floating point}. In particular,
in addition to the novel estimator we use $R^{(t)}$ (another bidiagonal
matrix formed by LSQR) instead of $D^{(t)}$ (as is done in LSQR),
which is more robust to inexact arithmetic. Using $R^{(t)}$ is somewhat
expensive, and this is why it is not used in LSQR, but because of
the other estimator, our algorithm needs to use this estimator only
once. Third, the other estimate that our algorithm computes comes
with a vector that proves that our estimate is an upper bound on $\sigma_{\min}$
(the approximate singular vector). LSQR's estimate also provides an
upper bound, but it does not provide a certificate vector which materializes
the bound. Fourth, the stopping criteria of our algorithm are designed
specifically to discover that the estimate of $\sigma_{\min}$ is
accurate; these criteria cannot be used in general in LSQR because
they require knowledge of the forward error.

\section{\label{sec:Mathematical-Motivation}The Forward Error of Iterates
on a Random Right-hand Side}

In this section we introduce the idea that underlies our algorithm,
namely that when applying LSQR to a random right-hand sides, the forward
error in the iterates is rich in the direction of the smallest singular
vector.

Let $x^{\star}$ be a random vector and let $b=Ax^{\star}$. Paige
and Saunders explain in the LSQR paper~\cite[Section~7.1]{LSQR}
that in exact arithmetic, the LSQR iterates are identical to the iterates
of Conjugate Gradients~\cite[Section~11.3]{golub13} when applied
to the normal equations
\[
A^{T}Ax^{\star}=A^{T}b\;.
\]
 Therefore, the $k$th iterate $x^{\left(k\right)}$ satisfies
\[
x^{\left(k\right)}=\arg\min_{x\in\mathcal{K}_{k}(A^{T}A,A^{T}b)}\Vert x^{\star}-x\Vert_{A^{T}A}
\]
 where in the above $\mathcal{K}_{k}(A^{T}A,A^{T}b)$ is the $k$-dimensional
Krylov subspace and $\left\Vert \;\cdot\;\right\Vert _{A^{T}A}$ denotes
the norm $\left\Vert v\right\Vert _{A^{T}A}^{2}=v^{T}A^{T}Av$. We
can therefore express $x^{(k)}$ as
\begin{align*}
x^{\left(k\right)} & =p_{0}(A^{T}A)A^{T}b
\end{align*}
for some degree-$(k-1)$ polynomial $p_{0}$ and the error as
\begin{align*}
x^{\star}-x^{\left(k\right)} & =x^{\star}-p_{0}(A^{T}A)A^{T}b\\
 & =x^{\star}-p_{0}(A^{T}A)A^{T}Ax^{\star}\\
 & =p(A^{T}A)x^{\star}\;,
\end{align*}
for a degree-$k$ polynomial $p(z)=1-zp_{0}(z)$.

Let us denote the thin SVD of $A$ and the eigendecomposition of its
Gram matrix by
\begin{eqnarray*}
A & = & U\Sigma V^{T}\\
A^{T}A & = & V\Sigma^{2}V^{T}\;,
\end{eqnarray*}
 where $\Sigma^{2}=\text{diag}(\sigma_{1}^{2},\sigma_{2}^{2},\dotsc,\sigma_{n}^{2})$
for some $\sigma_{1}\geq\sigma_{2}\geq\cdots\geq\sigma_{n}=\sigma_{\min}\geq0$.
Then
\begin{align*}
\min_{x^{(k)}}\Vert x^{\star}-x^{\left(k\right)}\Vert_{A^{T}A}^{2} & =\min_{p}\Vert p(A^{T}A)x^{\star}\Vert_{A^{T}A}^{2}\\
 & =\min_{p}\Vert p(V\Sigma^{2}V^{T})x^{\star}\Vert_{A^{T}A}^{2}\\
 & =\min_{p}\Vert Vp(\Sigma^{2})V^{T}x^{\star}\Vert_{A^{T}A}^{2}\\
 & =\min_{p}\sum_{i=1}^{n}\alpha_{i}^{2}\sigma_{i}^{2}p^{2}(\sigma_{i}^{2})\;,
\end{align*}
where $\alpha_{i}=(V^{T}x^{\star})_{i}$ and the minimization is
over all polynomials $p(z)$ of degree $k$ such that $p(0)=1$.

Now consider the $x^{(k)}$s, the forward errors $d^{(k)}=x^{\star}-x^{(k)}$,
and the Rayleigh quotients
\[
\frac{\|Ad^{(k)}\|}{\|d^{(k)}\|}\;.
\]
In subsequent sections we experimentally show that these quantities
tend to $\sigma_{\min}$ (and then diverge). To understand why this
happens, let us express these quotients as
\begin{eqnarray}
\frac{\|Ad^{(k)}\|^{2}}{\|d^{(k)}\|^{2}} & = & \frac{\|A(x^{\star}-x^{(k)})\|^{2}}{\|x^{\star}-x^{(k)}\|_{2}^{2}}\nonumber \\
 & = & \frac{\|AVp(\Sigma^{2})V^{T}x^{\star}\|^{2}}{\|Vp(\Sigma^{2})V^{T}x^{\star}\|^{2}}\nonumber \\
 & = & \frac{\|U\Sigma p(\Sigma^{2})V^{T}x^{\star}\|^{2}}{\|Vp(\Sigma^{2})V^{T}x^{\star}\|^{2}}\nonumber \\
 & = & \frac{\sum_{i=1}^{n}\alpha_{i}^{2}\sigma_{i}^{2}p^{2}(\sigma_{i}^{2})}{\sum_{j=1}^{n}\alpha_{j}^{2}p^{2}(\sigma_{j}^{2})}\label{eq:ratio-of-sums}\\
 & = & \sum_{i=1}^{n}\sigma_{i}^{2}\frac{\alpha_{i}^{2}p^{2}(\sigma_{i}^{2})}{\sum_{j=1}^{n}\alpha_{j}^{2}p^{2}(\sigma_{j}^{2})}\;.\nonumber
\end{eqnarray}

LSQR chooses $p$ among all degree-$k$ polynomials with $p(0)=1$
so as to minimize the numerator $\sum_{i=1}^{n}\alpha_{i}^{2}\sigma_{i}^{2}p^{2}(\sigma_{i}^{2})$
in~(\ref{eq:ratio-of-sums}), which is equal to the square of the
$A^{T}A$ norm of the forward error $d^{(k)}$. If $x^{\star}$ is
a random vector with independent standard independent elements, then
the $\alpha_{i}$s are also normal and independent, so we expect all
the $\alpha_{i}$s to be roughly similar in magnitude; we do not expect
very large or very small $\alpha_{i}$s. The term in the numerator
involving $p^{2}(\sigma_{i}^{2})$ is multiplied by $\alpha_{i}^{2}\sigma_{i}^{2}$.
If $\sigma_{n}<\sigma_{n-1}$, we expect the value of $p^{2}$ at
$\sigma_{1}^{2},\ldots,\sigma_{n-1}^{2}$ to be small relative to
$p^{2}(\sigma_{n}^{2})$, because the term involving $p^{2}(\sigma_{n}^{2})$
in the numerator of (\ref{eq:ratio-of-sums}) is scaled by $\sigma_{n}^{2}$
($\sigma_{n}=\sigma_{\min}$ is the smallest singular value), and
the $\alpha_{i}$s are all similar in magnitude. In other words, to
make the numerator of (\ref{eq:ratio-of-sums}) small, the LSQR polynomial
$p$ must assume small values at large $\sigma_{i}$s but can assume
relatively large values at small $\sigma_{i}$s.

If $p^{2}$ indeed assumes a small value at $\sigma_{1}^{2},\ldots,\sigma_{n-1}^{2}$
but a relatively large value at $\sigma_{n}^{2}$ and if the $\alpha_{i}^{2}$s
are all similar, then the expressions
\[
\frac{\alpha_{i}^{2}p^{2}(\sigma_{i}^{2})}{\sum_{j=1}^{n}\alpha_{j}^{2}p^{2}(\sigma_{j}^{2})}
\]
are much smaller than $1$ for $i=1,\ldots,n-1$ but close to $1$
for $i=n$ (these expressions sum to $1$). If this is the case, our
Rayleigh quotient approximates $\sigma_{n}$. The analysis generalizes
easily to multiple singular value at $\sigma_{\min}$.

We have not analyzed rigorously the quality of this approximation
as a function of the spectral gap (between $\sigma_{\min}$ and the
second smallest singular value) and of the random choice of $x^{\star}$.
Neverthess, this mathematical analysis motivates our algorithm, which
we describe in the next section.

\section{\label{sec:The-Algorithm}The Algorithm}

\begin{algorithm}
\begin{algorithmic}[1]

\small{

\STATE \textbf{Input: $A\in\mathbb{R}^{m\times n}$}

\STATE \textbf{Parameters (defaults in parentheses): $c_{1}$} ($8\epsilon_{\text{m}}$),
$c_{2}$ ($10^{-3}$), $c_{3}$ ($1/(64\cdot\epsilon_{\text{m}})$),
$c_{4}$ ($\sqrt{\epsilon_{\text{m}}}$), $c_{1}^{\prime}$ ($4\epsilon_{\text{m}}$),
and maxit.

\STATE

\STATE Estimate $\hat{\sigma}_{\max}=\sigma_{\max}$, along with
a certificate $\hat{v}_{\max}$, using power iteration.

\STATE $\hat{\sigma}_{\min}=\hat{\sigma}_{\max}$, $\hat{v}_{\min}=\hat{v}_{\max}$

\STATE Draw a random vector $\hat{x}\in\mathbb{R}^{n}$ with independent
normal entries

\STATE $\tau\gets\sqrt{2}\cdot\erf^{-1}(c_{2})/\Vert\hat{x}\Vert$

\STATE $x^{\star}\gets\hat{x}/\Vert\hat{x}\Vert$

\STATE $b\gets Ax^{\star}$

\STATE {*} $\beta^{(0)}\gets\Vert b\Vert$, $u^{(0)}\gets b/\beta^{(0)}$

\STATE {*} $v^{(0)}\gets Au^{(0)}$, $\alpha^{(0)}\gets\Vert v^{(0)}\Vert$,
$v^{(0)}\gets v^{(0)}/\alpha^{(0)}$

\STATE {*} $w^{(0)}\gets v^{(0)}$

\STATE {*} $x^{(0)}\gets0_{n\times1}$

\STATE {*} $\bar{\phi}^{(0)}\gets\beta^{(0)}$, $\bar{\rho}^{(0)}\gets\alpha$

\STATE $T\gets\text{maxit}$\qquad{} \COMMENT{limit for $T$; usually reduced later}

\FOR{$t=1,\dots,T$}

\STATE {*} $u^{(t)}\gets Av^{(t-1)}-\alpha u^{(t-1)}$

\STATE {*} $\beta^{(t)}\gets\Vert u^{(t)}\Vert$

\STATE {*} $u^{(t)}\gets u^{(t)}/\beta^{(t)}$

\STATE {*} $v^{(t)}\gets A^{T}u^{(t)}-\beta^{(t)}v^{(t)}$

\STATE {*} $\alpha^{(t)}\gets\Vert v^{(t)}\Vert$

\STATE {*} $v^{(t)}\gets v^{(t)}/\alpha^{(t)}$

\STATE {*} $\rho^{(t)}\gets\left\Vert \left(\begin{array}{cc}
\bar{\rho}^{(t-1)} & \beta^{(t)}\end{array}\right)\right\Vert $

\STATE {*} $c^{(t)}\gets\bar{\rho}^{(t-1)}/\rho^{(t)}$, $s^{(t)}\gets\beta^{(t)}/\rho^{(t)}$

\STATE {*} $\theta^{(t)}\gets s^{(t)}\alpha^{(t)}$, $\bar{\rho}^{(t)}\gets-c^{(t)}\alpha^{(t)}$

\STATE {*} $\phi^{(t)}\gets c^{(t)}\bar{\phi}^{(t-1)}$, $\bar{\phi}^{(t)}\gets s\bar{\phi}^{(t-1)}$

\STATE {*} $x^{(t)}\gets x^{(t-1)}+(\phi^{(t)}/\rho^{(t)})w^{(t-1)}$

\STATE {*} $w^{(t)}\gets v^{(t)}-(\theta^{(t)}/\rho^{(t)})w^{(t-1)}$

\STATE $R_{tt}^{(t)}\gets\rho^{(t)}$\qquad{} \COMMENT{Only diagonal and superdiagonal of $R$ are kept in memory}

\STATE \textbf{if $t>1$ set $R_{t-1,t}^{(t)}\gets\theta^{(t-1)}$}

\STATE $d^{(t)}\gets x^{\star}-x^{(t)}$

\STATE \textbf{if }$d^{(t)}=0$ \textbf{set $\hat{\sigma}_{\min}\gets\hat{\sigma}_{\max}$,
$\hat{v}_{\min}\gets\hat{v}_{\max}$ and break for}\\
 \qquad{} \COMMENT{Due to finite precision arithmetic, when $\kappa \neq 1$ getting $d^{(t)}=0$ is extremely unlikely}.

\IF{$\Vert A d^{(t)} \Vert \leq \hat{\sigma}_{\min} \Vert d^{(t)} \Vert$}

\STATE $\hat{\sigma}_{\min}\gets\Vert Ad^{(t)}\Vert/\Vert d^{(t)}\Vert,$
$\hat{v}_{\min}\gets d^{(t)}$

\ENDIF

\STATE  \textbf{if $\hat{\sigma}_{\min}/\hat{\sigma}_{\max}\leq c_{4}$
then }$c_{1}\gets c_{1}^{\prime}$

\IF{ $T=\text{maxit}$ and ($\frac{\Vert Ad^{(t)}\Vert}{\hat{\sigma}_{\max}\Vert x^{(t)}\Vert+\Vert b\Vert}\leq c_{1}$
\textbf{or }$\Vert d^{(t)}\Vert\leq\tau$ \textbf{or $\hat{\sigma}_{\max}/\hat{\sigma}_{\min}\geq c_{3}$})
}

\STATE\textbf{$T\gets\left\lceil 1.25t\right\rceil $}

\ENDIF

\ENDFOR

\STATE Estimate $\tilde{\sigma}_{\min}=\sigma_{\min}(R^{(T)})$,
using inverse power iteration

\STATE $\tilde{\sigma}_{\min}\gets\min(\tilde{\sigma}_{\min},\hat{\sigma}_{\min})$

\STATE \RETURN $\hat{\sigma}_{\max}/\hat{\sigma}_{\min}$, $\hat{\sigma}_{\max}$,$\hat{\sigma}_{\min}$,
$\hat{v}_{\max}$ and $\hat{v}_{\min}$, $\tilde{\sigma}_{\min}$

}

\end{algorithmic}

\caption{\label{alg:the-alg}The condition number estimation algorithm. Lines
marked with {*} are inherited from LSQR. Some quantities, such as
$\Vert\hat{x}\Vert$, $\Vert Ad^{(t)}\Vert$, $\Vert d^{(t)}\Vert$,
and $\Vert b\Vert$, should be computed once and reused; the full
expressions are shown multiple times for clarity.}
\end{algorithm}

This section describes our algorithm for estimating the condition
number of $A$. A detailed pseudo-code description appears in Algorithm~\ref{alg:the-alg}.

The algorithm starts by estimating $\sigma_{\max}$ and a corresponding
certificate vector using power iteration on $A^{T}A$. By a \emph{certificate
vector }for an estimate $\hat{\sigma}_{\max}$ we mean a vector $v$
such that $\|Av\|/\|v\|=\hat{\sigma}_{\max}$ (and similarly for estimates
of $\sigma_{\min}$). We perform enough iterations to estimate $\sigma_{\max}$
to within 10\% of accuracy with probability at least $1-10^{-12}$.
Using a bound due to Klein and Lu~\cite[ Section 4.4]{KleinLu1996}\footnote{Note that the statement of Lemma~6 in~\cite{KleinLu1996} is incorrect;
the proof shows the correct bound. Also, the discussion that follows
the proof of the lemma repeats the error in the statement of the lemma.}, we find that given a relative error parameter $\epsilon$ and a
failure probability parameter $\delta$, if we perform
\[
\left\lceil \frac{1}{\epsilon}\left(\ln\left(2n\right)^{2}+\ln\left(\frac{1}{\epsilon\delta^{2}}\right)\right)\right\rceil
\]
iterations, the relative error in our approximation is less than $\epsilon$
with probability at least $1-\delta$. For the parameters $\epsilon=10^{-1}$
and $\delta=10^{-12}$, 1004 iterations suffice even for matrices
with up to $10^{9}$ columns. For $\epsilon=1/3$ and $\delta=10^{-12}$,
only 298 iterations suffice for matrices with up to $10^{9}$ columns.
(The accuracy of the $\sigma_{\max}$ estimate in the power method
is typically much higher than predicted by this bound, but the additional
accuracy depends on the gap between the largest and second-largest
singular values; the bound that we use makes no assumption on the
gap.) It is possible to use other methods to find $\hat{\sigma}_{\max}$
and $v$ (e.g. Lanczos). We chose power iteration because when we
developed the algorithm existing theoretical results allow us to write
an expression for the number of iterations as a function of $\epsilon$
and $\delta$, independent of the gap. Such strong bounds that guarantee
high accuracy of $\sigma_{\max}$ are desirable since the stopping
criteria we use for the estimation of $\sigma_{\min}$ depend on $\hat{\sigma}_{\max}$.
Furthermore, for one singular pair, we empirically observed that the
power method is typically faster than Lanczos and other iterative
methods.

The main phase of the algorithm uses a modified LSQR iteration~\cite[pages 50--51]{paige82,LSQR}
designed to estimate $\sigma_{\min}$ and to produce a corresponding
certificate vector. LSQR is a method for solving least squares problems
$\min\|Ax-b\|_{2}$. At its core, LSQR uses the bidiagonlization procedure
of Golub and Kahan to form iterates $u^{(0)},u^{(1)},\dots\in\mathbb{R}^{m}$,
$v^{(0)},v^{(1)},\dots\in\mathbb{R}^{n},\alpha^{(0)},\alpha^{(1)},\dots\in\mathbb{R}$
and $\beta^{(0)},\beta^{(1)},\dots\in\mathbb{R}$ such that
\begin{eqnarray*}
U^{(t+1)}(\beta^{(1)}e_{1}) & = & b\\
AV^{(t)} & = & U^{(t+1)}B^{(t)}\\
A^{T}U^{(t+1)} & = & V^{(t)}(B^{(t)})^{T}+\alpha^{(t+1)}v^{(t+1)}e_{t+1}^{T}
\end{eqnarray*}
 where $U^{(t)}=\left[u^{(0)},u^{(1)},\cdots,u^{(t)}\right]$ and
$V_{k}=\left[v^{(0)},v^{(1)},\cdots,v^{(t)}\right]$. LSQR uses these
to find at iteration $t$ the optimal minimizer inside the Krylov
subspace ${\cal K}_{t}(A^{T}A,A^{T}b)=span\{A^{T}b,(A^{T}A)A^{T}b,\dots,(A^{T}A)^{t-1}A^{T}b\}.$
In addition, LSQR forms additional bidiagonal matrices $R^{(0)},R^{(1)},\dots$
and estimates $\bar{\phi}^{(0)},\bar{\phi}^{(1)},\dots$ of the norm
of the residual $r^{(t)}=b-Ax^{(t)}$. See~\cite{LSQR} for details.
We also remark that the Golub and Kahan bidiagonlization is closely
related to the Lanczos tridiagonlization of the augmented matrix and
of the Gram matrix, and that this connection is exploited by our algorithm.

The main idea in our algorithm is to use a random right-hand side
so the forward error will tend to concentrate in the direction corresponding
to $\sigma_{\min}$. To be able to access the forward error, we also
generate the right-hand side in a way that gives us access to it.
Specifically, our algorithm first generates a uniformly-distributed
random vector $x^{\star}$ on the unit sphere by first generating
a vector $\hat{x}$ with normally-distributed independent random components,
and setting $x^{\star}=\hat{x}/\|\hat{x}\|$. The algorithm multiplies
it by $A$ to produce a nearly consistent right-hand side $b$. While
in exact arithmetic we would have had a consistent right hand side
($b=Ax^{\star}$), due to finite precision it is only nearly consistent:
$\|b-Ax^{\star}\|\leq O(\epsilon_{\text{m}})\cdot\|A\|$.

Our algorithm modifies LSQR by adding a few steps to each iteration.
At the end of each (standard) LSQR iteration, we have an updated approximate
solution $x^{(t)}$ and an estimate of $\|r^{(t)}\|=\|Ax^{(t)}-b\|$,
denoted by $\bar{\phi}^{(t)}$. The equality of $\bar{\phi}^{(t)}$
and $\|r^{(t)}\|$ depends on orthogonality of the Lanczos vectors,
which lose orthogonality in floating point arithmetic as the algorithm
progresses, so in practice $\bar{\phi}^{(t)}$ and $\|r^{(t)}\|$
may differ significantly. Our algorithm also computes the forward
error $d^{(t)}=x^{\star}-x^{(t)}$ and $\|d^{(t)}\|$. We have
\begin{eqnarray*}
\left\Vert Ad^{(t)}\right\Vert  & = & \left\Vert A\left(x^{\star}-x^{(t)}\right)\right\Vert \\
 & = & \left\Vert b-Ax^{(t)}\right\Vert \\
 & = & \left\Vert r^{(t)}\right\Vert \;,
\end{eqnarray*}
which in exact arithmetic equals $\bar{\phi}^{(t)}$, but to improve
the robustness of the algorithm we compute $\|Ad^{(t)}\|$ explicitly.
(In our numerical experiments we have found $\bar{\phi}^{(t)}$ to
be an accurate estimate, but we prefer to avoid any reliance on the
orthogonality of the Lanczos vectors in our algorithm.) We also compute
$\|x^{(t)}\|$.

Next, the algorithm computes the ratio $\|Ad^{(t)}\|/\|d^{(t)}\|,$
which like any Rayleigh quotient is an upper bound on $\sigma_{\min}$.
If this ratio is the smallest we have seen so far, the algorithm treats
it as an estimate of $\sigma_{\min}$ and stores both the ratio and
the certificate $d^{(t)}$. When the algorithm terminates, it outputs
the best ratio it has found and the corresponding certificate.

We use three stopping criteria. The first criterion, which is the
one used by the standard LSQR algorithm for consistent systems, stops
the algorithm when the \emph{normwise backward error}~\cite[Section~7.1]{Higham02}
drops below a threshold,
\begin{equation}
\frac{\left\Vert r^{(t)}\right\Vert }{\hat{\sigma}_{\max}\left\Vert x^{(t)}\right\Vert +\left\Vert b\right\Vert }\leq c_{1}\,,\label{eq:stop1}
\end{equation}
where $\hat{\sigma}_{\max}$ is our estimate of $\Vert A\Vert$ and
$c_{1}$ is a parameter that is set by default to $8\epsilon_{\text{m}}$.
It has been observed experimentally~\cite{CPT09} that for consistent
systems, as long as $c_{1}$ is of the order of magnitude of $\epsilon_{\text{m}}$
or greater,  this criterion will be eventually met in spite of the
loss of orthogonality in the biorthogonalization process; however,
the left-hand side of \eqref{eq:stop1} does not seem to decrease
much below the value required to satisfy the inequality~\cite{CPT09}.

In many cases our second stopping criterion will stop the algorithm
well before the residual is that small. This second condition is
\begin{equation}
\left\Vert d^{(t)}\right\Vert \leq\frac{\sqrt{2}\cdot\erf^{-1}(c_{2})}{\Vert\hat{x}\Vert}\;,\label{eq:stop2}
\end{equation}
where $\erf^{-1}$ is the inverse error function, computed using a
numerical approximation, and $c_{2}$ is a parameter that is set by
default to $10^{-3}$, which corresponds to $\erf^{-1}(c_{2})\approx8.9\times10^{-4}$.
This criterion is not used by LSQR; it requires knowledge of the error,
which LSQR does not have in general. We explain this stopping criterion,
and how the choice of $c_{2}$ affects the algorithm, in the next
section.

The third stopping criterion is
\begin{equation}
\frac{\hat{\sigma}_{\max}}{\hat{\sigma}_{\min}}\geq c_{3}\,,\label{eq:stop3}
\end{equation}
where $c_{3}$ is a parameter that is set by default to $1/(64\cdot\epsilon_{\text{m}})$.
In other words, at this threshold we consider the matrix to be numerically
rank deficient and we do not attempt to estimate the exact condition
number. Standard LSQR uses a different condition number estimate in
a similar stopping criterion, for regularization~\cite{LSQR}.

To achieve good accuracy even for matrices that are terribly ill conditioned
(condition number close to $1/\epsilon_{\text{m}}$), the stopping
criteria are refined in two additional ways:
\begin{enumerate}
\item If at some point we have
\[
\frac{\hat{\sigma}_{\max}}{\hat{\sigma}_{\min}}\geq c_{4}\,,
\]
where $c_{4}$ is a parameter that is set by default to $\sqrt{\epsilon_{\text{m}}}$,
we set $c_{1}$ (residual-based stopping threshold) to $c_{1}^{\prime}$,
which is set by default to $4\epsilon_{\text{m}}$. The rationale
is that when the matrix is ill conditioned, the residual based test
is less accurate so we lower threshold to compensate.
\item Even when the method detects convergence using one of its three criteria
(small residual, small error, and numerical rank deficiency), it keeps
iterating. The number of extra iterations is one quarter of the number
performed until convergence was detected. This rule is a heuristic
that tries to improve the accuracy of the condition number estimate.
The cost of this heuristic is obviously limited and it can be turned
off by the user.
\end{enumerate}

In addition, our algorithm also stores the matrix $R^{(t)}$, one
of two bidiagonal matrices that LSQR incrementally constructs but
normally discards. In exact arithmetic, the singular values of $R^{(t)}$
converge to the singular values of $A$~\cite{cullum85,druskin91}.
Once the algorithm terminates, we compute $\tilde{\sigma}_{\min}\approx\sigma_{\min}(R^{(t)})$;
if it is smaller than the best $\|Ad^{(t)}\|/\|d^{(t)}\|$ estimate,
we output both estimates. One estimate ($\tilde{\sigma}_{\min}$)
is tighter, but it comes with no certificate vector; the other is
looser, but comes with a certificate. Generating the certificate for
the Lanczos estimate requires storing the Lanczos vectors or repeating
the iterations. The former is obviously too expensive, while the latter
might be redundant if the user is satisfied with the certificate already
computed. We note that the singular values of $R^{(t)}$ tend to converge
first to large singular values of $A$, so this estimate is not very
likely to be better than $\|Ad^{(t)}\|/\|d^{(t)}\|$.

Storing $R^{(t)}$ and estimating $\sigma_{\min}(R^{(t)})$ is relatively
inexpensive since $R^{(t)}$ is bidiagonal. We estimate it by running
inverse iteration on $R^{(t)}$, again performing enough iterations
to get to within 10\% of accuracy with very high probability ($1-10^{-12}$;
this is a parameter in the code). The cost of a single inverse power
iteration is only $O(t)$. Since we use power-iteration, the error
in $\tilde{\sigma}_{\min}$ is one sided: it is always the case that
$\tilde{\sigma}_{\min}\geq\sigma_{\min}(R^{(t)})$ (because it is
generated by a Rayleigh quotient). Also, $\sigma_{\min}(R^{(t)})\geq\sigma_{\min}$.
Therefore, $\tilde{\sigma}_{\min}$ is also an upper bound on $\sigma_{\min}$,
not just an estimate. We remark that $\tilde{\sigma}_{\min}$ is guaranteed
to be an upper bound on $\sigma_{\min}$ only in exact arithmetic.
In floating point, it is possible for $\tilde{\sigma}_{\min}$ to
drop below $\sigma_{\min}$.

\section{\label{sec:analysis-small-err}Analysis of the Small-Error Stopping
Criterion}

We now explain the second stopping criteria \eqref{eq:stop2}. Suppose
that the smallest singular value of $A$ is simple and that it is
well separated from the next smallest singular values. Let $x^{\star}=\sum_{i=1}^{n}\alpha_{i}v_{i}$
be the initial vector represented in the basis of the right singular
vectors of $A$. As LSQR progresses towards finding $x^{\star}$ it
will tend initially to resolve components in the direction of the
largest singular vectors. Since the $v_{n}$ direction is not present
in $x^{(0)}=0$ we expect it to be not present during the initial
iterations, i.e., $(v_{n})^{T}x^{(t)}\approx0$. This implies that
we expect for the initial iterations to have $\vert(v_{n})^{T}(x^{\star}-x^{(t)})\vert\approx\vert\alpha_{n}\vert$
so $\|x^{\star}-x^{(t)}\|\geq\alpha_{n}$. Now, at some point in the
iteration, the solution $x^{(t)}$ will be roughly $x^{(t)}\approx\sum_{i=1}^{n-1}\alpha_{i}v_{i}$,
i.e., the error remains mostly in the direction of the smallest singular
subspace, but the $v_{n}$ direction is not present at all. At that
point, $\|x^{\star}-x^{(t)}\|\approx\vert\alpha_{n}\vert$. LSQR will
now start to resolve that error at least partially and the norm of
the error will decrease below $\vert\alpha_{n}\vert$. If we stop
the iteration when $\|x^{\star}-x^{(t)}\|\gg\vert\alpha_{n}\vert$,
the error is unlikely to be a good estimate of a small singular vector.
If we stop when $\|x^{\star}-x^{(t)}\|\leq\vert\alpha_{n}\vert$ we
will likely have a good estimate of a small singular value. Ideally,
we want to stop immediately when $\|x^{\star}-x^{(t)}\|$ drops below
$\vert\alpha_{n}\vert$. Stopping later (when the error is much smaller
than $\vert\alpha_{n}\vert$) does not do any harm, since we report
the best Rayleigh quotient seen, but it does not improve the estimate
by much.

The second stopping criteria is designed so that the test passes only
if the condition $\|x^{\star}-x^{(t)}\|\leq\vert\alpha_{n}\vert$
holds with high probability. It is based on the following proposition.
\begin{prop}
Suppose that $x^{\star}=\hat{x}/\|\hat{x}\|$ where $\hat{x}$ is
a vector with normally-distributed independent random components,
and assume $x^{\star}=\sum_{i=1}^{n}\alpha_{i}v_{i}$ be $x^{\star}$
represented in the basis of the right singular vectors of $A$. If
\[
\|x^{\star}-y\|\leq\frac{\sqrt{2}\cdot\erf^{-1}(\delta)}{\Vert\hat{x}\Vert}
\]
 then $\|x^{\star}-y\|\leq\vert\alpha_{n}\vert$ holds with probability
of at least $1-\delta$ for $0<\delta<1$.
\end{prop}

\begin{proof}
Let $\hat{y}=\Vert\hat{x}\Vert y$ and $\hat{\alpha}_{n}=\Vert\hat{x}\Vert\alpha_{n}$.
Obviously, $\|x^{\star}-y\|\leq\vert\alpha_{n}\vert$ if and only
if $\|\hat{x}-\hat{y}\|\leq\vert\hat{\alpha}_{n}\vert$. Therefore,
if we find a value $\tau$ such that $\Pr(\vert\hat{\alpha}_{n}\vert\geq\tau)\geq1-\delta$
then the condition $\|\hat{x}-\hat{y}\|\leq\tau$ implies that $\|\hat{x}-\hat{y}\|\leq\vert\hat{\alpha}_{n}\vert$,
and hence $\|x^{\star}-y\|\leq\vert\alpha_{n}\vert$, with a probability
of at least $1-\delta$. Since $\Vert v_{n}\Vert=1$ we have $\hat{\alpha}_{n}\sim N(0,1)$.
Therefore,
\[
\Pr(\vert\hat{\alpha}_{n}\vert\geq\sqrt{2}\cdot\erf^{-1}(\delta))=1-\delta\,,
\]
so $\tau=\sqrt{2}\cdot\erf^{-1}(\delta)$ can be used. The condition
$\|\hat{x}-\hat{y}\|\leq\tau$ is now equivalent to $\|x^{\star}-y\|\leq\tau/\Vert\hat{x}\Vert$
which is exactly the criteria in the proposition statement.
\end{proof}
The choice of $c_{2}$ in the algorithm determines our confidence
that $\|x^{\star}-x^{(t)}\|$ dropped below $\vert\alpha_{n}\vert$.
We use a default value of $10^{-3}$, which implies a small probability
of failure, but not a tiny one. The user can, of course, change the
value if he or she needs higher confidence (in either case the error
is one sided). Another approach is to use a much larger $c_{2}$ and
to repeat the algorithm $\ell$ times (possibly in a single run, exploiting
matrix-matrix multiplies), say $\ell=3$. The probability that we
succeed in at least one run is at least $1-c_{2}^{\ell}$. For $c_{2}=10^{-2}$,
say, setting $\ell=3$ or so should suffice. In our experience it
is better to make $c_{2}$ smaller than to set $\ell>1$, but we did
not do a formal analysis of this issue.

If the smallest singular value is multiple (associated with a singular
subspace of dimension $k>1$), our situation is even better, because
we can stop when $x^{\star}-x\approx\sum_{i=n-k+1}^{n}\alpha_{i}v_{i}$,
when $\|x^{\star}-x\|\approx\sqrt{\alpha_{n-k+1}^{2}+\cdots+\alpha_{n}^{2}}$,
which is even more likely to be larger than our stopping criterion.

When the smallest singular value is not well separated, the stopping
criterion is still sound but the Rayleigh quotient estimate we obtain
is not as accurate, because in such cases $x^{\star}-x$ tends to
be a linear combination of singular vectors corresponding to several
singular values. These singular values are all small, but they are
not exactly the same, thereby pulling the Rayleigh quotient up a bit.
\textbf{\textcolor{red}{}}

\section{\label{sec:numerical}Numerical Investigations}

This section illustrate and explore the behavior of our algorithm
using numerical experiments.

\subsection{\label{sec:Rationale}Experiments on Synthetic Matrices }

We begin by exploring the behavior of our method on synthetic $1000$-by-$400$
full-rank real matrices with prescribed singular values and random
singular vectors.

The first matrix has 10 singular values at $10^{-8}$, 300 that are
distributed logarithmically between $10^{-3}$ and $10^{-2}$, and
90 at $1$. The large gap between the smallest singular value and
the next larger one makes the problem relatively easy for our algorithm.
\begin{figure}
\begin{centering}
\includegraphics[width=0.6\textwidth]{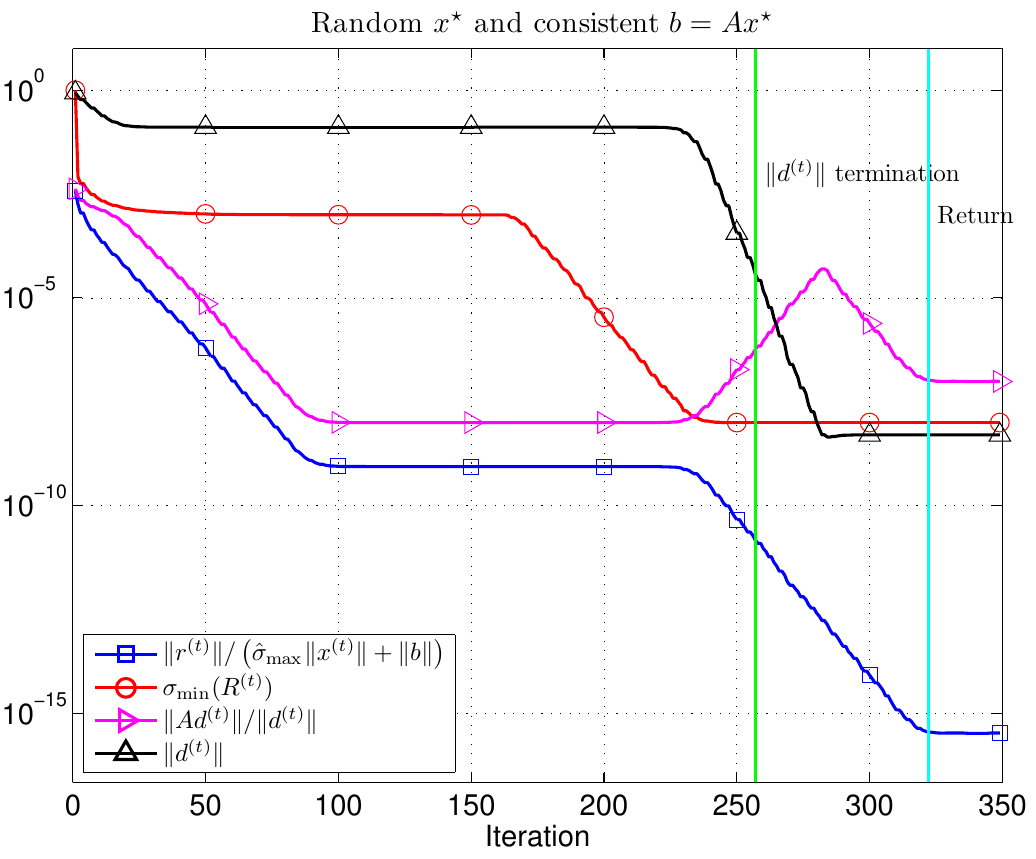}
\par\end{centering}
\caption{\label{fig:illustration_our}The behavior of the method on one matrix
described in the text. The vertical lines indicate the iterations
at which the stopping criterion was triggered (green), and the iteration
at which the algorithm would have returned (after performing 25\%
extra iterations; we ran the algorithm longer in this case).}
\end{figure}
 Figure~\ref{fig:illustration_our} shows the behavior of our method
on one matrix generated with this spectrum. In the first 90 iterations
or so the backward error (left-hand side in \eqref{eq:stop1}) diminishes
logarithmically; the norm of $x^{(t)}-x^{\star}$ drops a bit but
then stops dropping much. These two effects cause our estimate $\hat{\sigma}_{\min}$
to also diminish roughly logarithmically. The Lanczos estimate $\sigma_{\min}(R^{(t)})$
drops a bit initially but stagnates from iteration 30 or so. What
is happening up to iteration 90 or so is that Golub-Kahan bidiagonalization
resolves the singular values in the $10^{-3}$-to-$10^{-2}$ cluster,
while LSQR removes much of the projection of the corresponding singular
vectors from the residual and from the error. Around iteration 160
Lanczos has resolved enough of the spectrum in the $10^{-3}$-to-$10^{-2}$
cluster and the smallest singular value of $R^{(t)}$ starts moving
toward $10^{-8}$. At that point, most of the remaining error consists
of singular vectors corresponding to the singular value $10^{-8}$,
which causes our estimate to be accurate (to within 9 decimal digits!).
The norm of $x^{(t)}-x^{\star}$ is still large, because the error
contains a significant component in the subspace associated with the
$10^{-8}$ singular values. At this point, around iteration 225, when
the topmost black curve starts dropping again, LSQR starts to resolve
the error in this subspace, the backward error starts decreasing again,
the norm of $x^{(t)}-x^{\star}$ starts decreasing, which causes stopping
criterion~\eqref{eq:stop2} to be met about 30 iterations later.

\begin{figure}
\begin{centering}
\includegraphics[width=0.6\textwidth]{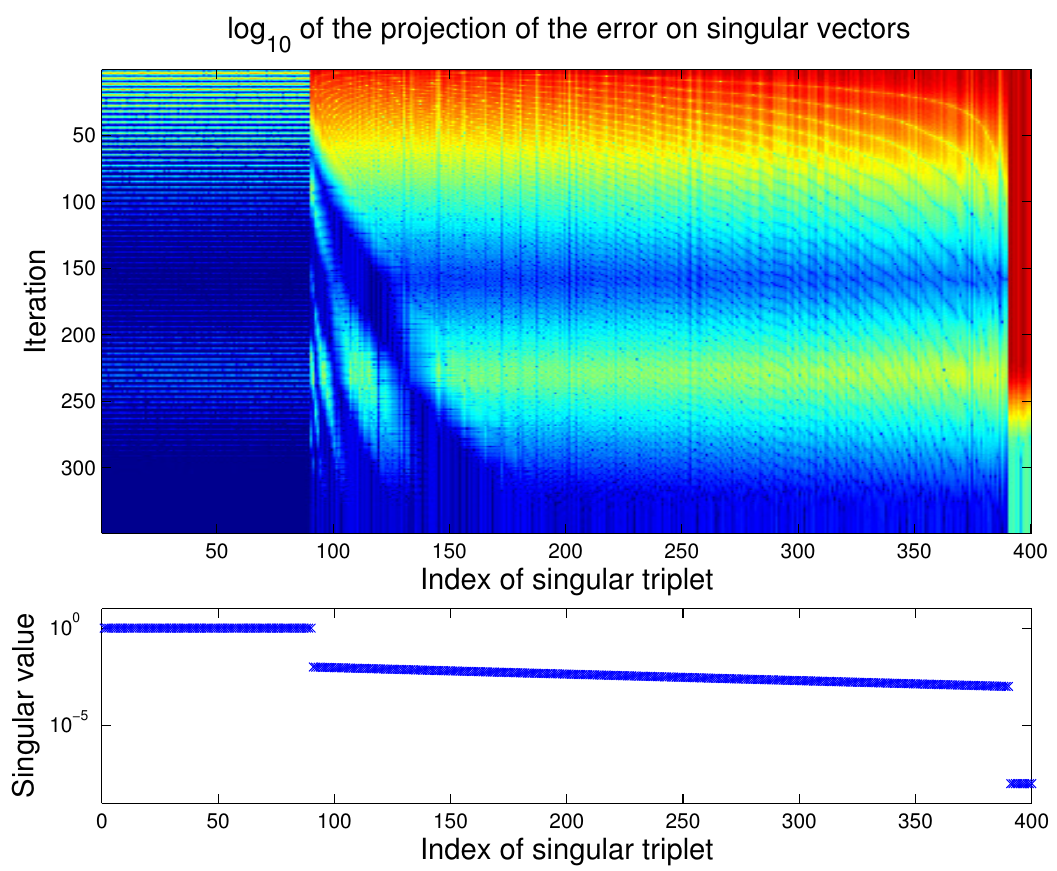}
\par\end{centering}
\caption{\label{fig:illustration_our_proj}The projection of the forward error
on the right singular vectors in the experiment shown in Figure~\ref{fig:illustration_our}.
The bottom plot shows the singular values of $A$ and the top image
shows the projections. Blue represents values near $\epsilon_{\text{m}}$;
green represents values near $\sqrt{\epsilon_{\text{m}}}$; and red
represents values near $1$.}
\end{figure}
Figure~\ref{fig:illustration_our_proj} visualizes the behavior described
in the previous paragraph by plotting the projection of the forward
error on the right singular vectors. We see that up to around iteration
50, the error associated with singular spaces associated with singular
values other than $1$ remains very large. From that point on to about
iteration 150, the error in subspaces corresponding to values between
$10^{-3}$ and $10^{-2}$ is resolved, but the error associated with
the singular value $10^{-8}$ is still very large. This is a point
where our method finds the smallest singular value and its certificate
(the error). As LSQR starts to resolve the error in the $10^{-8}$
singular subspace, the errors in the $10^{-3}$ to $10^{-2}$ subspaces
grows (perhaps due to loss of orthogonality) but they are reduced
again later.

The method yielded similar results when the smallest singular value
was moved down to $10^{-13}$, with convergence after about 440 iterations
and an estimate that is correct to within 5 decimal digits. The value
$\sigma_{\min}=10^{-13}$ is about the lower limit for which the $\|d^{(t)}\|$
stopping criterion is useful.

\begin{figure}
\begin{centering}
\includegraphics[width=0.6\textwidth]{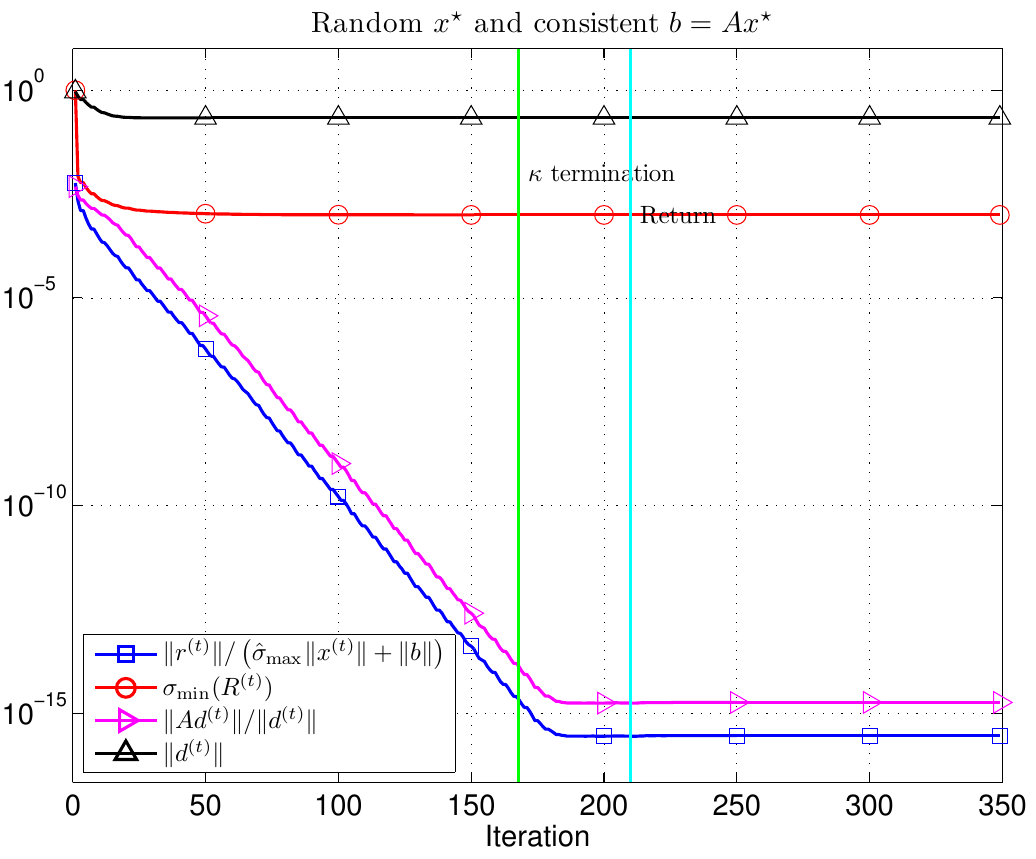}
\par\end{centering}
\caption{\label{fig:illustration_rankdef_our}The behavior of the method on
a numerically rank deficient matrix described in the text. Termination
here occurred due to criterion~\eqref{eq:stop3}.}
\end{figure}

When $A$ is rank deficient there are infinitely many solutions to
the system $Ax=b$. The solution $x^{\star}$, which was generated
randomly, has no special property that distinguishes it from other
solutions (like minimum norm), so no least-squares solver can recover
$x^{\star}$. Therefore, it is unlikely $\|d^{(t)}\|$ will become
small enough to cause our method to stop. In this case, the method
stops because the residual eventually becomes very small (close to
$\epsilon_{\text{m}}$) or because the estimated condition number
becomes too big (stopping condition~\eqref{eq:stop3}). In Figure~\ref{fig:illustration_rankdef_our}
we illustrate the behavior of the algorithm on a rank deficient matrix;
the matrix has 10 singular values that are $10^{-16}$ (numerical
zeros), 300 distributed logarithmically between $10^{-3}$ and $10^{-2}$,
and 90 at $1$.
The algorithm stopped because the condition number got too big; however
a few iterations later it would have stopped because $\Vert r^{(t)}\Vert$
became too small (stopping condition~\eqref{eq:stop1}). Our estimate
is not accurate (around $1.8\times10^{-15}$, a relative error of
$18$), but it still indicates to the user that the matrix is numerically
rank deficient. Stopping criteria \eqref{eq:stop1} and \eqref{eq:stop3},
and the relatively-inaccurate estimates they yield, are triggered
only when the matrix is close to rank deficiency (condition number
of about $10^{14}$ or larger).

Can
\[
\|\arg\min_{x}\|Ax-b\|_{2}\|_{2}^{-1}
\]
estimate $\sigma_{\min}$ if we take $b$ to be random unit vector?
If $A$ has fewer columns than rows or is rank deficient, then with
high probability a random $b$ is not in the column space of $A$.
On some matrices, the minimizer has a norm that is larger than the
norm of $b$ by about a factor of $\sigma_{\min}^{-1}$. However,
unless there is a large gap between the smallest singular values and
the rest, the minimizer has a smaller norm and the norms ratio fails
to accurately estimate $\sigma_{\min}^{-1}$.
\begin{figure}
\begin{centering}
\includegraphics[width=0.6\textwidth]{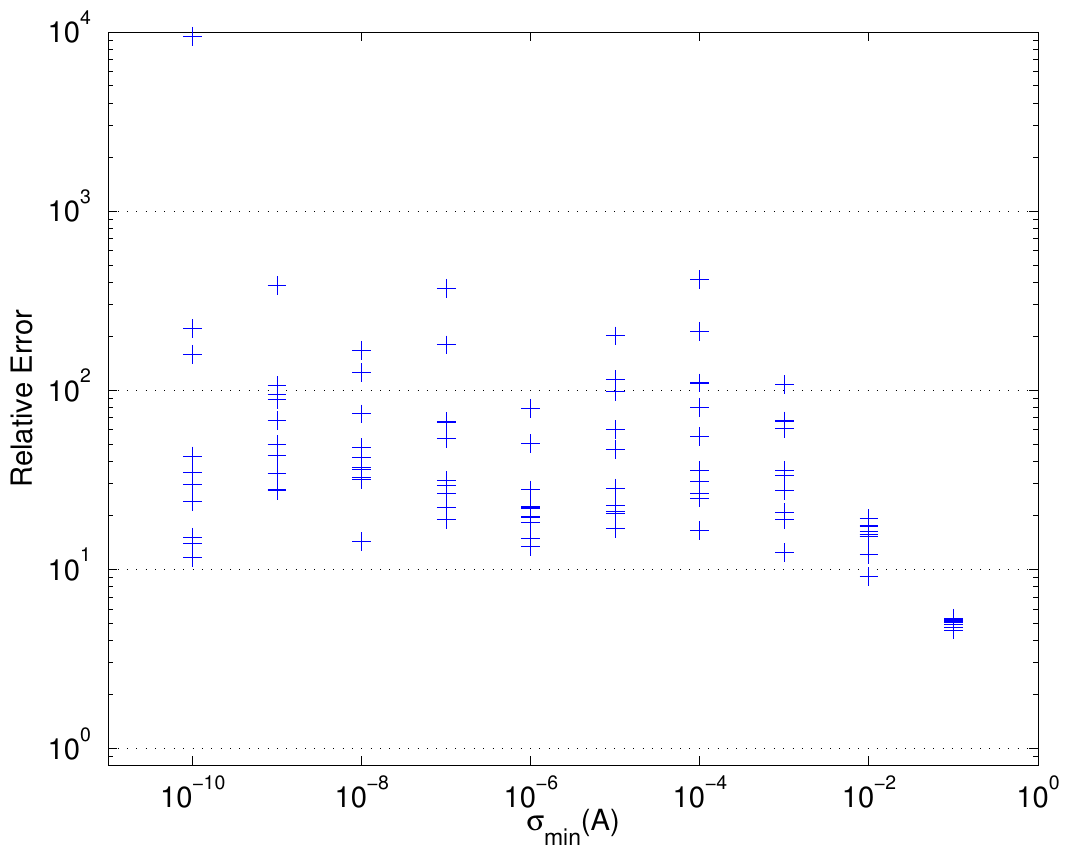}
\par\end{centering}
\caption{\label{fig:random_x_norm_of_b}Estimating $\sigma_{\min}$ by $\|A^{+}b\|^{-1}$
for a unit-length random $b$ with normal independent components.}
\end{figure}
 Figure~\ref{fig:random_x_norm_of_b} shows the relative errors in
this estimate for matrices whose singular values are distributed linearly
between $\sigma_{\min}$ and $1$. The errors are huge, sometimes
by more than 3 orders of magnitude. This is not a particularly useful
method. This method amounts to one half of inverse iteration on $A^{T}A$,
so it is not surprising that it is not accurate; performing more iterations
would make the method more reliable, but at the cost of applying the
pseudoinverse many times.

This estimator is clearly biased (the estimate is always larger than
$\sigma_{\min}$); so is any fixed number of steps of inverse iteration.
Kenney et al.~\cite{Kenney98} derive an unbiased estimator of this
type for the Frobenius-norm condition number. To the best of our knowledge,
this is not possible in the Euclidean norm.

Other distributions of the singular values lead to more accurate estimates
in this method. But will this method work with an approximate minimizer
produced by an iterative method? The following experiment suggests
that the answer is no, at least for LSQR. The matrix used in the experiment
has 50 singular values distributed logarithmically between $10^{-10}$
and $10^{-9}$, 50 more distributed logarithmically between $10^{-1}$
and $1$, and the rest are all $1$.
\begin{figure}
\begin{centering}
\includegraphics[width=0.6\textwidth]{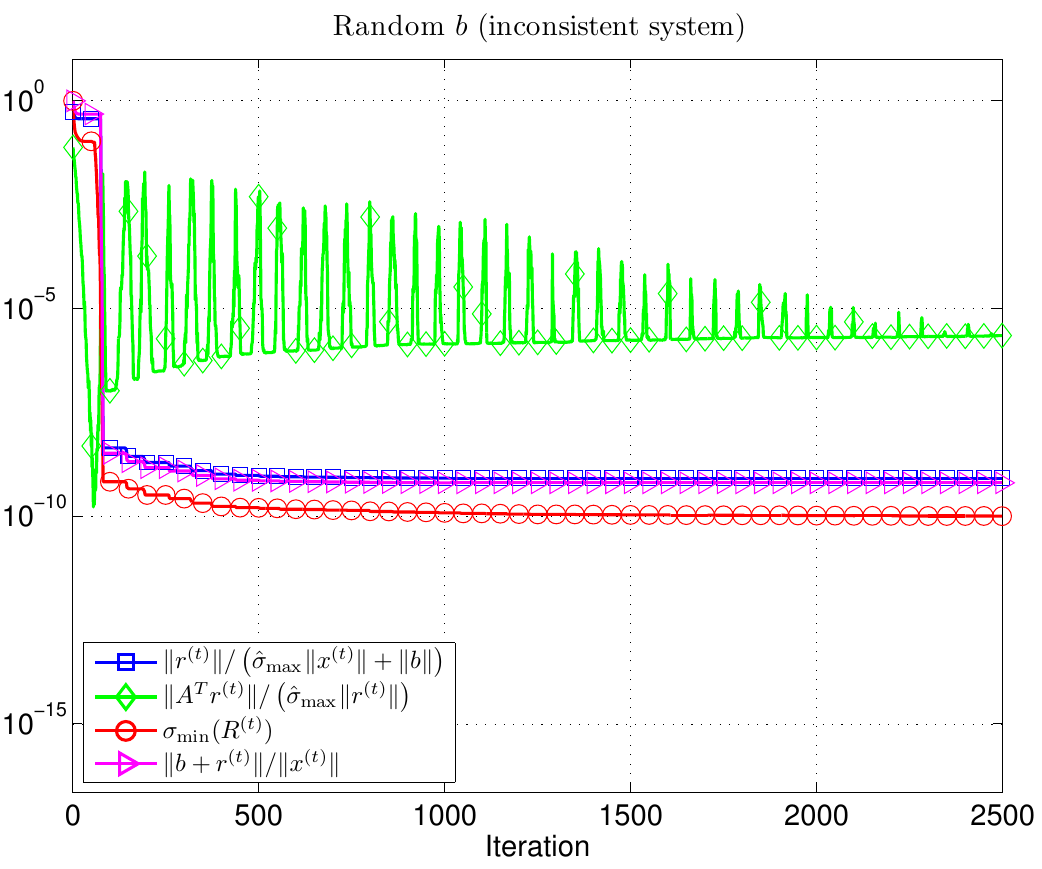}
\par\end{centering}
\caption{\label{fig:no_Atr_convergence}Running LSQR with an inconsistent right-hand-side
$b$.}
\end{figure}
The results, presented in Figure~\ref{fig:no_Atr_convergence}, indicate
that there is no good way to decide when to terminate LSQR when used
in this way to estimate $\sigma_{\min}$. We obviously cannot rely
on the residual approaching $\epsilon_{\text{m}}$, because the problem
is inconsistent. The original LSQR paper~\cite{LSQR} suggests another
stopping condition,
\[
\frac{\left\Vert A^{T}r^{(t)}\right\Vert }{\hat{\sigma}_{\max}\left\Vert r^{(t)}\right\Vert }\leq c\;,
\]
but our experiment shows that this ratio may fail to get close to
$\epsilon_{\text{m}}$. In our experiment the best local minimum is
around $10^{-10}$, six orders of magnitude larger than $\epsilon_{\text{m}}$.
Moreover, at that local minimum, around iteration 60, the estimate
$\|b+r^{(t)}\|/\|x^{(t)}\|$ is still near $1$, very far from $\sigma_{\min}$.
The Lanczos estimate $\sigma_{\min}(R^{(t)})$ is also very inaccurate
at that time. There does not appear to be a good way to decide when
to stop the iterations and to report the best estimate seen so far.

\begin{figure}
\begin{centering}
\includegraphics[width=0.6\columnwidth]{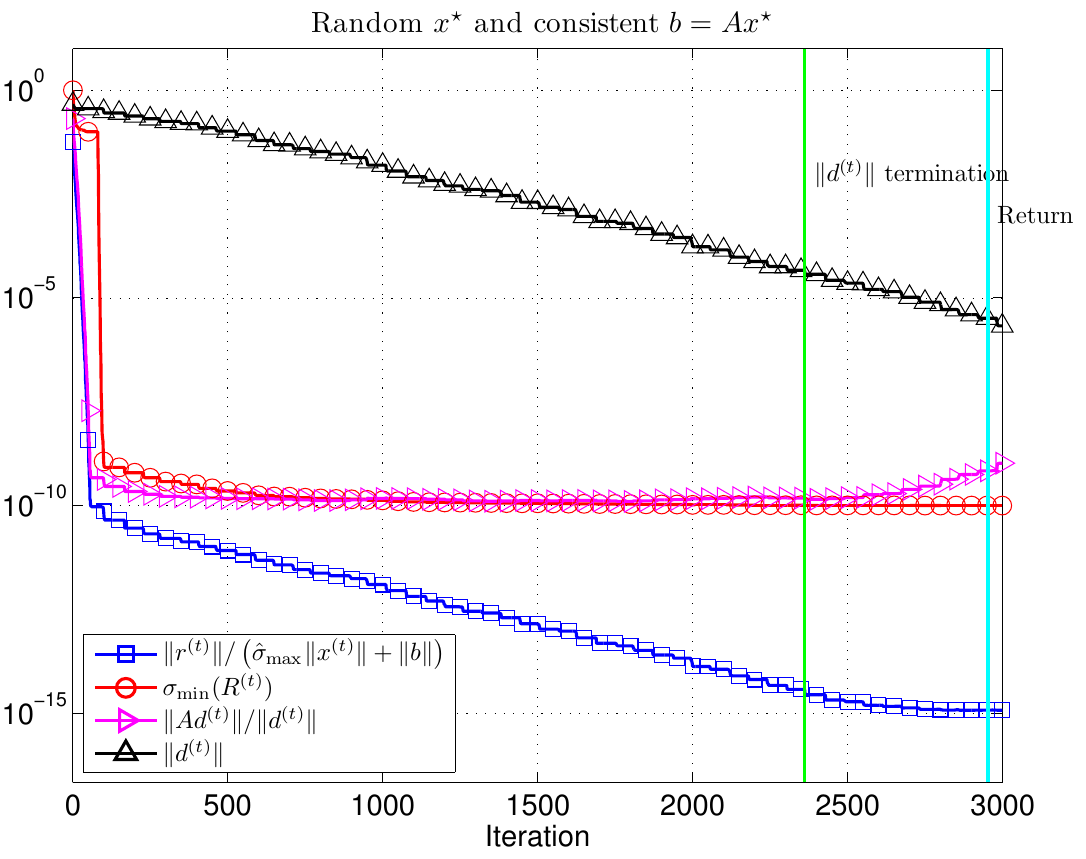}
\par\end{centering}
\caption{\label{fig:no_Atr_convergence_our}Running our method on the same
matrix as in Figure~\ref{fig:no_Atr_convergence}.}
\end{figure}

On matrices with this singular value distribution, our method detects
convergence after 2400-2500 iterations, returning a certified estimate
of $\sigma_{\min}$ that is accurate to within 15\textendash 40\%
(the accuracy of the Lanczos estimate is better, with relative errors
smaller than 10\%). Figure~\ref{fig:no_Atr_convergence_our} shows
a typical run. The number of iterations is large, but the stopping
criteria are robust.

\begin{figure}
\begin{centering}
\includegraphics[width=0.6\textwidth]{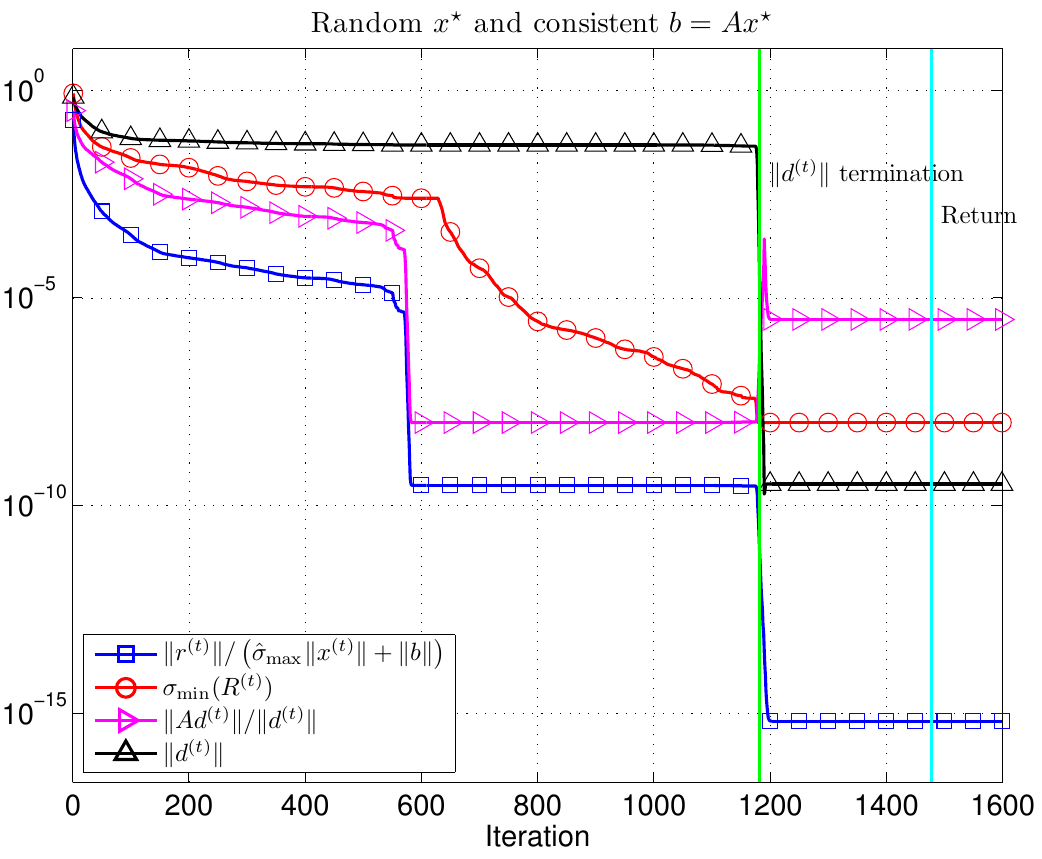}
\par\end{centering}
\caption{\label{fig:linear_1e-8}Singular values are distributed linearly from
$10^{-8}$ up to $1$.}
\end{figure}

In Figure~\ref{fig:linear_1e-8} all the singular values are distributed
linearly from $10^{-8}$ up to $1$ (the dimensions of all matrices
in this subsection are again $1000$-by-$400$). Convergence is fairly
slow. The gap between $\sigma_{\min}=\sigma_{400}=10^{-8}$ and $\sigma_{399}$
is relatively large, around $\frac{1}{400}$, so $\sigma_{\min}$
is computed accurately.
 When many singular values are distributed logarithmically or nearly
so, convergence is very slow and the small relative gap between $\sigma_{\min}$
and the next-larger singular values causes the method to return a
less accurate estimate.
\begin{figure}
\begin{centering}
\includegraphics[width=0.6\textwidth]{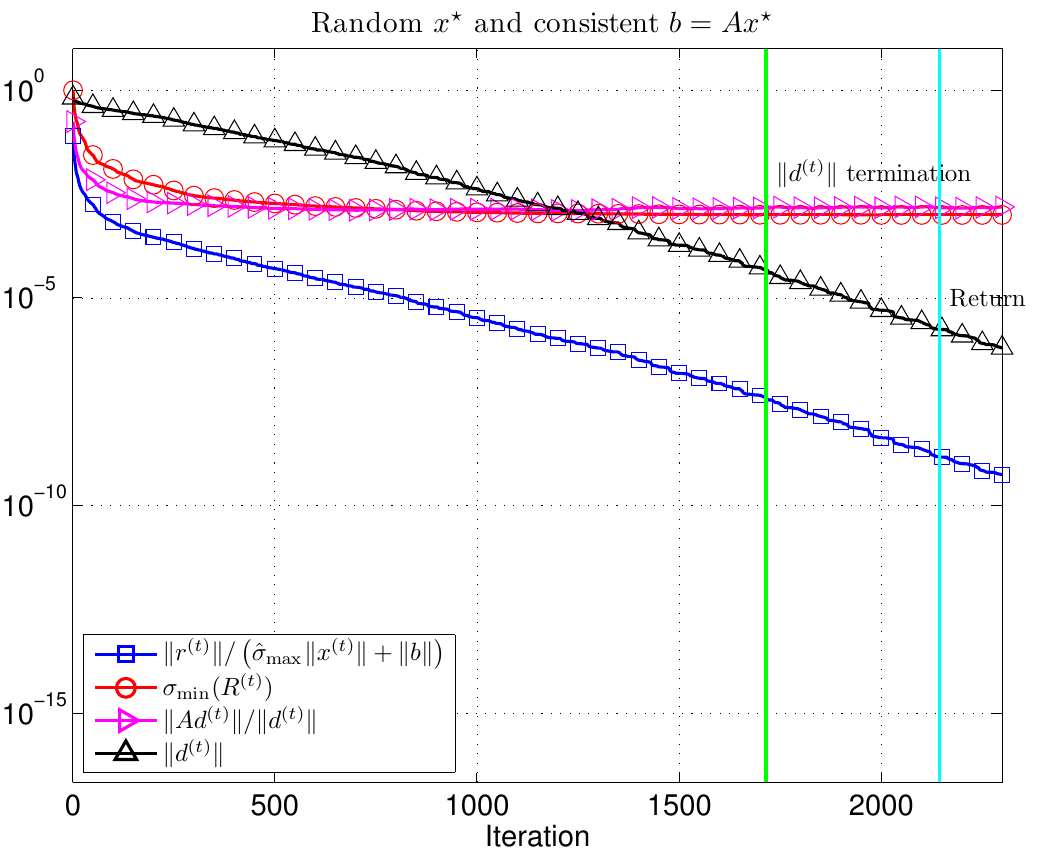}
\par\end{centering}
\caption{\label{fig:log_1e-3}200 singular values are distributed logarithmically
from $10^{-3}$ up to $1$, and the rest are equal to $1$.}
\end{figure}
Figure~\ref{fig:log_1e-3} plots the convergence when 200 singular
values are distributed logarithmically between $10^{-3}$ and $1$
and the rest are at $1$. We do not see a period of stagnation during
which the error is a good estimate of $v_{\min}$. The certified estimate
is only accurate to within 31\% and the Lanczos estimate to within
10\% (much worse than when the smallest singular value is well separated
from the rest).

\begin{figure}
\begin{centering}
\includegraphics[width=0.6\textwidth]{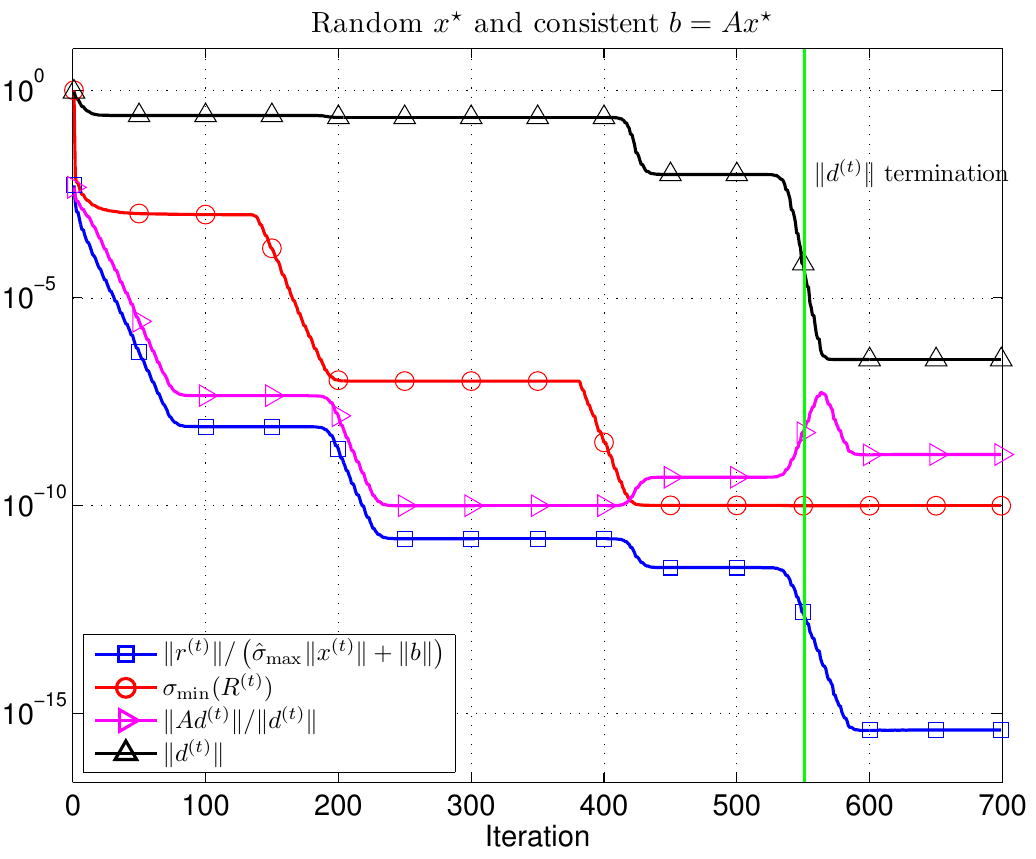}
\par\end{centering}
\caption{\label{fig:multiple_stagnations}200 singular values are distributed
logarithmically from $10^{-3}$ up to $1$, and the rest are equal
to $1$.}
\end{figure}

LSQR might have several periods of stagnation. This happens when the
spectrum contains several well-separated clusters. Figure~\ref{fig:multiple_stagnations}
plots the convergence when the matrix has a multiple singular value
at $10^{-10}$, a multiple singular value at $10^{-7}$ (both with
multiplicity 10), 300 singular values that are distributed logarithmically
between $10^{-3}$ and $10^{-2}$, and the rest are at $1$. We see
multiple stagnation periods of both the residual, the error, the Lanczos
estimate, and our certified estimate.

\subsection{Experiments on Large Structured Random Matrices}

The next set of experiments was performed on sparse matrices that
motivated this project. These matrices have exactly three nonzeros
per column, where the location of the nonzeros is random and uniform
and their values are $+1$ or $-1$ with equal independent probabilities.
This type of matrices arises when simulating the evolution of random
2-dimensional complexes in various stochastic models~\cite{ALLM12}.
Such $m$-by-$n$ matrices tend to be well conditioned when $n<0.9m$
and rank deficient when $n>0.95m$.
\begin{figure}
\noindent \begin{centering}
\includegraphics[width=0.6\textwidth]{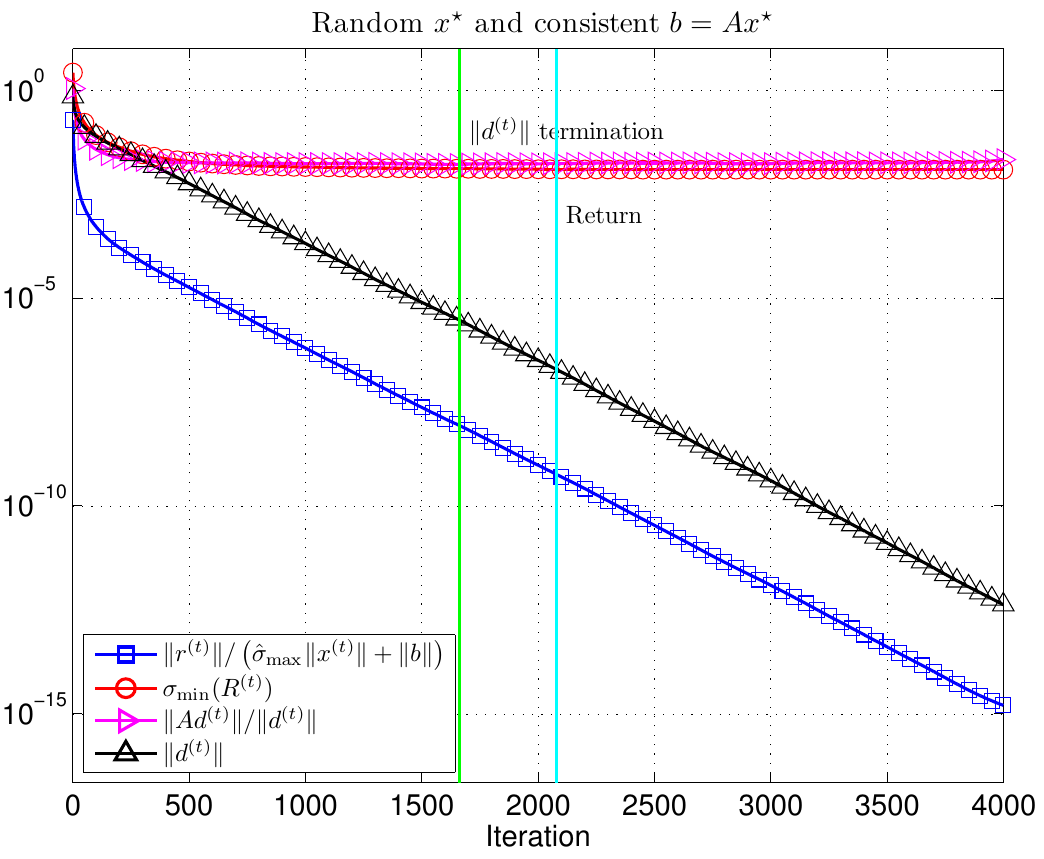}
\par\end{centering}
\bigskip{}
\noindent \begin{centering}
\includegraphics[width=0.6\textwidth]{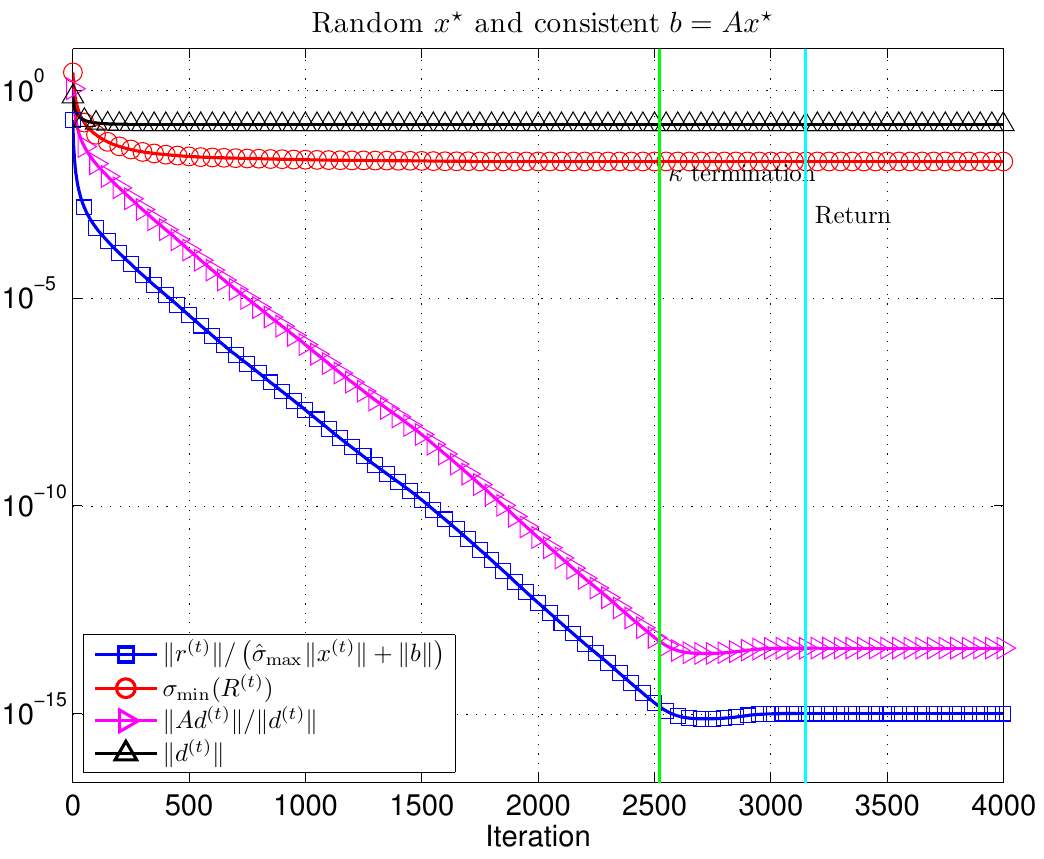}
\par\end{centering}
\caption{\label{fig:irad_m=00003D100000}Random matrices with 3 nonzeros per
row. On the top we see the convergence on a $100{,}000$-by-$90{,}000$
matrix and on the bottom convergence on a $100{,}000$-by-$95{,}000$
matrix.}
\end{figure}
Figure~\ref{fig:irad_m=00003D100000} shows that the method converges
quite quickly even on large matrices in both the well conditioned
and the rank deficient cases. On smaller matrices of this type we
were able to assess the accuracy of the method. For $m=1000$, on
$n=900$ the algorithm yielded a relative error of 22\% (the Lanczos
estimate was off by 78\%), and on $n=450$ the algorithm yielded a
relative error of 41\% (the Lanczos estimate was off by 18\%). Problems
of this type of size $m=1{,}000{,}000$ required similar number of
iterations and were easily solved on a laptop. It is worth noting
that due to the random structure of the non-zero pattern, it is likely
that factorization based condition number estimators will be very
slow when applied to this type of matrices.

\subsection{\label{subsec:Large-Scale-Experiments}Experiments on Many Real-World
Matrices}

\begin{figure}
\begin{centering}
\includegraphics[width=0.47\columnwidth]{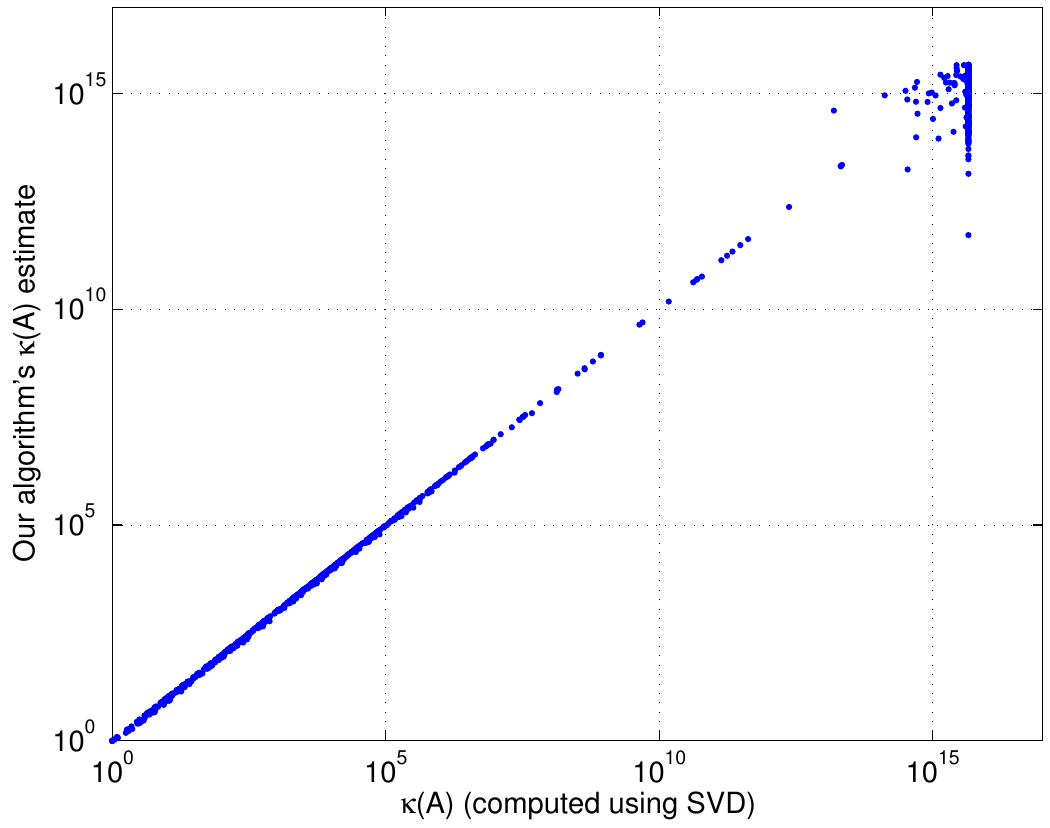}\hfill{}\includegraphics[width=0.47\columnwidth]{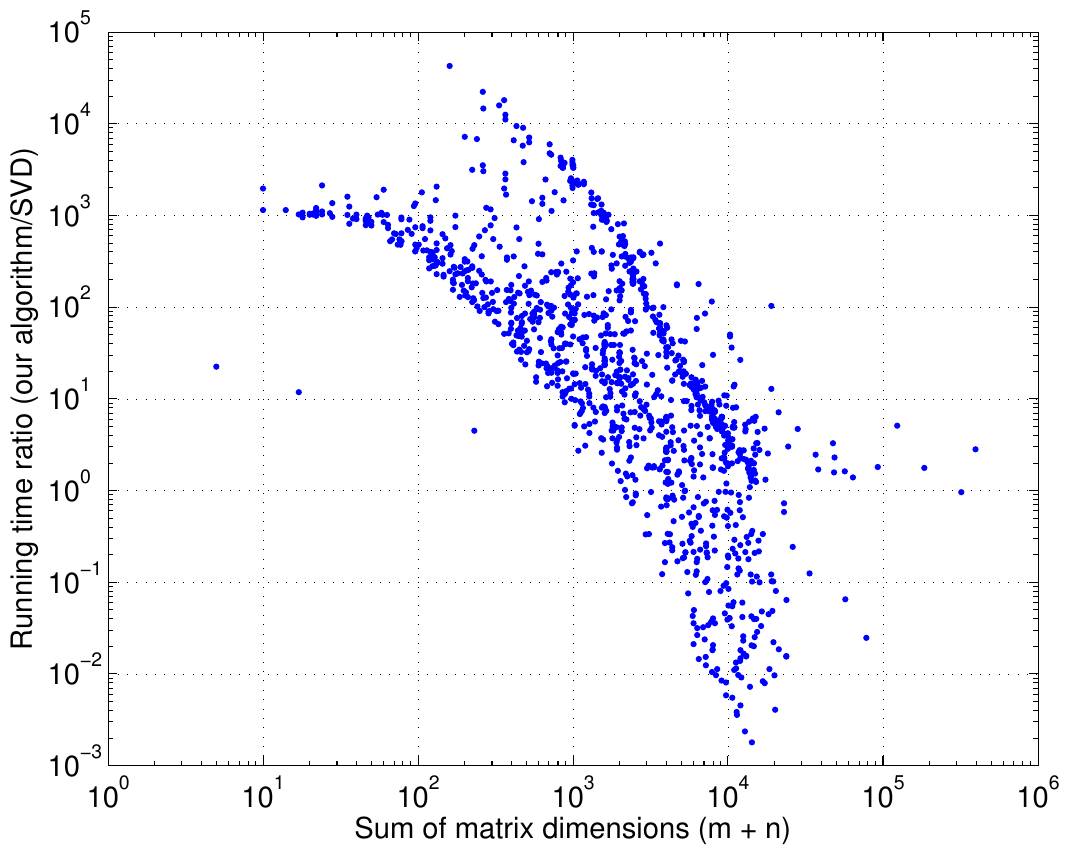}
\par\end{centering}
\caption{\label{fig:scatterplots}Correlation between the actual condition
number and the one that we obtain from our algorithm when it converges
(left) and the computational cost of our method (implemented in Matlab),
relative to the dense SVD, as a function of matrix size (right).}
\end{figure}
\begin{figure}
\begin{centering}
\includegraphics[width=0.47\columnwidth]{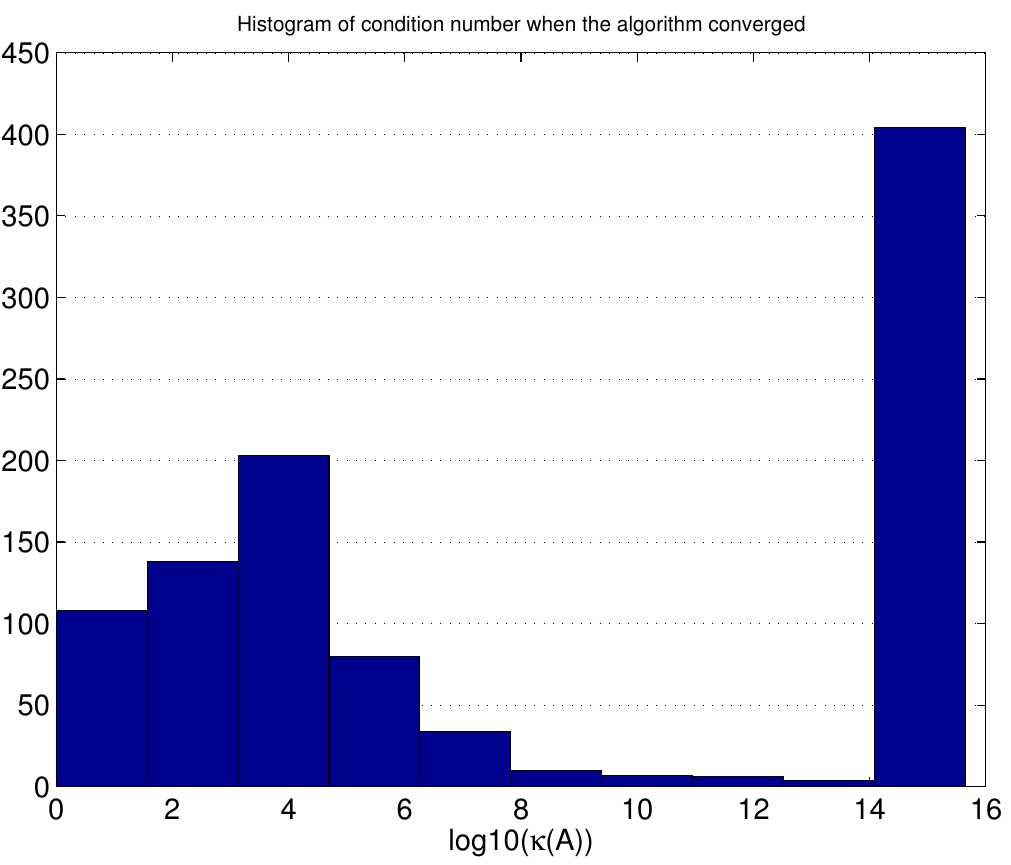}\hfill{}\includegraphics[width=0.47\columnwidth]{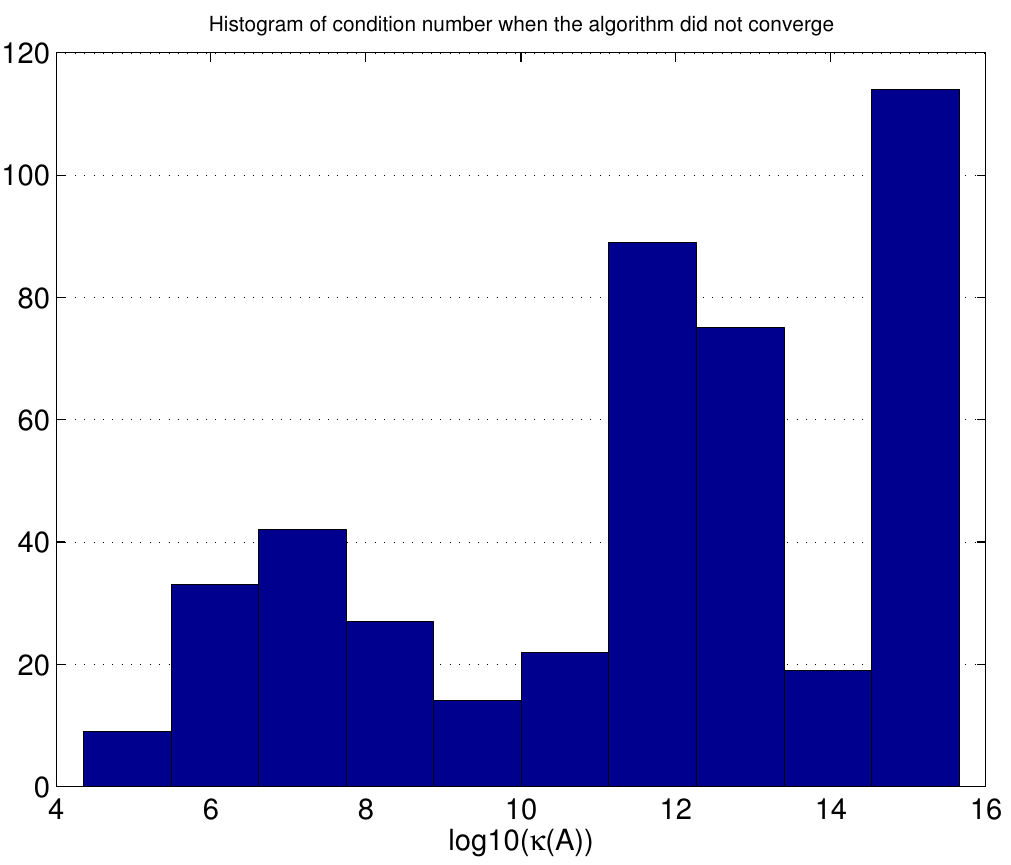}
\par\end{centering}
\caption{\label{fig:histograms}Histogram of the condition number of matrices
on which our algorithm converged (left) and did not converge (right).}
\end{figure}

We ran both a dense SVD and our method on all the matrices from Tim
Davis's sparse matrix collection~\cite{UFDavis11} for which $mn^{2}<256\times10^{9}$.
Our method converged in 100,000 iterations or less on 1024 out of
the 1468 matrices in this category.

Out of the 1468 matrices, 404 had condition number $1/(64\cdot\epsilon_{\text{m}})\approx7\times10^{13}$
or larger. Our method converged on 278 out of them, delivering condition
number estimates of $5\times10^{11}$ or larger. In other words, on
all the matrices that were close to rank deficiency, our method detected
that the condition number is large, but in some cases it underestimated
the actual condition number.

On matrices with condition number smaller than $64/\epsilon_{\text{m}}$,
our method always estimated the condition number to within a relative
error of 24\% or less. Figure~\ref{fig:scatterplots} shows a scatter
plot of the condition number reported by our method when it converged
vs. the condition number reported by Matlab's dense SVD.

We ran the method again on some of the matrices on which it failed
to converge in 100,000 iterations, allowing the method to run longer.
It converged in all cases. For example, on \texttt{nos1}, the method
detected convergence after 169,791 iterations. The $\sigma_{\min}$
estimate it returned was actually from iteration 90,173, meaning that
at iteration 100,000 it actually converged, but the algorithm was
not yet able to detect convergence. We note that \texttt{nos1} is
a square matrix of dimension 237; the method can be slow even on small
matrices. Figure~\ref{fig:histograms} presents histograms of the
dense-SVD condition-numbers of matrices on which our method converged
or did not converge. The histograms show that our method can fail
to converge (in 100,000 iterations) on both well condition and ill
conditioned matrices. It also appears that full-rank ill conditioned
matrices are more likely to cause our method difficulties than either
well conditioned or rank deficient matrices.

We remark the one clearly sees the effect of inexact arithmetic in
these experiments. Since LSQR on an $m$-by-$n$ matrix is equivalent
to CG on an $n$-by-$n$ matrix, with exact arithmetic our method's
iteration count should never exceed the dimension. In practice, we
see that the number of iterations may well exceed the dimension (e.g.,
the \texttt{nos1} matrix). This is because of loss of orthogonality
due to floating point rounding.

The running time of our method obviously varies a lot and is not easy
to characterize. But on large matrices it is often much faster than
a dense SVD, as shown in Figure~\ref{fig:scatterplots}.

\subsection{\label{subsec:Experiments-on-Large}Experiments on Large Real-World
Matrices}

We ran the algorithm on a few very large matrices from the same matrix
collection. As on the smaller matrices, the method sometimes converged,
although sometimes it exceeded the maximal number of iterations (up
to 1,000,000). Table~\ref{tab:very-large} shows a sample of the
the statistics of successful runs; they indicate that when the singular
spectrum is clustered, the method works well even on very large matrices.

\begin{table}
\begin{centering}
\begin{tabular}{lrr|rrr|rr|rr}
 &  &  & {\footnotesize{}condest} &  &  & {\footnotesize{}irlba} &  & {\footnotesize{}PRIMME} & \tabularnewline
 & {\footnotesize{}$m$} & {\footnotesize{}$n$} & {\footnotesize{}time (s)} & {\footnotesize{}it. (\#mv)} & {\footnotesize{}$\kappa$ (est.)} & {\footnotesize{}\#mv} & {\footnotesize{}$\kappa$ (est.)} & {\footnotesize{}\#mv} & {\footnotesize{}$\kappa$ (est.)}\tabularnewline
\cline{2-10}
{\footnotesize{}rajat10} & {\footnotesize{}30202} & {\footnotesize{}30202} & {\footnotesize{}3.7} & {\footnotesize{}3017 (6634)} & {\footnotesize{}1.1e+03} & {\footnotesize{}43016} & {\footnotesize{}1.3e+03} & {\footnotesize{}5248} & {\footnotesize{}1.3e+03}\tabularnewline
{\footnotesize{}flower\_7\_4} & {\footnotesize{}67593} & {\footnotesize{}27693} & {\footnotesize{}0.6} & {\footnotesize{}129 (858)} & {\footnotesize{}1.5e+01} & {\footnotesize{}104} & {\footnotesize{}2.0e+01} & {\footnotesize{}402} & {\footnotesize{}1.9e+01}\tabularnewline
{\footnotesize{}flower\_8\_4} & {\footnotesize{}125361} & {\footnotesize{}55081} & {\footnotesize{}3.2} & {\footnotesize{}648 (1896)} & {\footnotesize{}2.8e+13} & {\footnotesize{}4584} & {\footnotesize{}2.4e+16} & {\footnotesize{}3556} & {\footnotesize{}3.7e+00}\tabularnewline
{\footnotesize{}\uline{wheel\_601 }} & {\footnotesize{}902103} & {\footnotesize{}723605} & {\footnotesize{}358} & {\footnotesize{}5260} & {\footnotesize{}1.3e+14} & \multicolumn{2}{c|}{{\footnotesize{}FAIL}} & \multicolumn{2}{c}{{\footnotesize{}FAIL}}\tabularnewline
{\footnotesize{}\uline{Franz11}} & {\footnotesize{}47104} & {\footnotesize{}30144} & {\footnotesize{}0.3} & {\footnotesize{}73 (746)} & {\footnotesize{}$\infty$} & {\footnotesize{}712} & {\footnotesize{}3.3e+16} & {\footnotesize{}528} & {\footnotesize{}4.9e+00}\tabularnewline
{\footnotesize{}lp\_ken\_18} & {\footnotesize{}154699} & {\footnotesize{}105127} & {\footnotesize{}10} & {\footnotesize{}2234 (5068)} & {\footnotesize{}2.5e+14} & {\footnotesize{}169096} & {\footnotesize{}5.7e+16} & {\footnotesize{}24832} & {\footnotesize{}6.5e+02}\tabularnewline
{\footnotesize{}lp\_pds\_20} & {\footnotesize{}108175} & {\footnotesize{}33874} & {\footnotesize{}3.3} & {\footnotesize{}848 (2296)} & {\footnotesize{}1.4e+14} & {\footnotesize{}7304} & {\footnotesize{}2.8e+16} & {\footnotesize{}6178} & {\footnotesize{}6.8e+00}\tabularnewline
{\footnotesize{}cage15} & {\footnotesize{}5154859} & {\footnotesize{}5154859} & {\footnotesize{}162} & {\footnotesize{}73 (746)} & {\footnotesize{}1.1e+01} & {\footnotesize{}72} & {\footnotesize{}1.2e+01} & {\footnotesize{}92} & {\footnotesize{}1.2e+01}\tabularnewline
{\footnotesize{}atmosmodl} & {\footnotesize{}1489752} & {\footnotesize{}1489752} & {\footnotesize{}663} & {\footnotesize{}6699 (13998)} & {\footnotesize{}8.5e+02} & {\footnotesize{}662536} & {\footnotesize{}1.1e+03} & {\footnotesize{}15546} & {\footnotesize{}1.1e+03}\tabularnewline
{\footnotesize{}Rucci1} & {\footnotesize{}1977885} & {\footnotesize{}109900} & {\footnotesize{}2194} & {\footnotesize{}20576 (41752)} & {\footnotesize{}6.5e+03} & {\footnotesize{}1728264} & {\footnotesize{}6.8e+03} & {\footnotesize{}28604} & {\footnotesize{}6.8e+03}\tabularnewline
{\footnotesize{}LargeRegFile} & {\footnotesize{}2111154} & {\footnotesize{}801374} & {\footnotesize{}146} & {\footnotesize{}1304 (3208)} & {\footnotesize{}1.0e+04} & {\footnotesize{}21032} & {\footnotesize{}1.1e+04} & {\footnotesize{}4886} & {\footnotesize{}1.1e+04}\tabularnewline
{\footnotesize{}sls} & {\footnotesize{}1748122} & {\footnotesize{}62729} & {\footnotesize{}119} & {\footnotesize{}848 (2296)} & {\footnotesize{}9.7e+02} & {\footnotesize{}8296} & {\footnotesize{}1.3e+03} & {\footnotesize{}3218} & {\footnotesize{}1.3e+03}\tabularnewline
\end{tabular}
\par\end{centering}
\caption{\label{tab:very-large}Large real-world matrices whose condition number
was successfully computed by our method. Matrices with an underline
indicate a structurally rank deficient matrix, or a matrix indicated
in the matrix collection as numerically rank deficient. }
\end{table}

We also attempted to approximate the condition number using two singular
values solvers: irlba~\cite{BR05} and PRIMME~\cite{PRIMME_SVDS}.
Both are general purpose singular values solvers, so in order to make
a fair comparison with our algorithm (which is designed for condition
number estimation) we judiciously chose the stopping criteria. From
PRIMME, we use the library's option to supply a special convergence
criteria, and stop the method once the residual has dropped below
$10\%$ the current estimate of $\sigma_{\min}$. For irlba, we set
the code's convergence tolerance to $10^{-15}$ when our code detected
rank deficiency, and to $1/\kappa_{\text{est }}$ otherwise, where
$\kappa_{\text{est }}$ is the condition number estimate obtained
using our method (of course, this stopping criteria cannot be used
in practice, and we use it only for the sake of empirical evaluation).
For both algorithms, we restrict the number of iterations to 100,000.
Since the different libraries are implemented in different languages
(irlba and our algorithm are implemented in MATLAB, while PRIMME does
most of the work in C), we measure only the number of matrix-vector
products. For the matrices reported in Table~\ref{tab:very-large},
unless the matrix is very well conditioned, our algorithm clearly
outperforms irlba. As for PRIMME, when the condition number is not
too large, PRIMME tends to require less matrix-vector products. However,
PRIMME fails to detect very ill-conditioned matrices.

\section{Summary}

We have presented an adaptation of LSQR to the estimation of the condition
number of matrices.

Our method is yet another tool in the spectral condition-number estimation
toolbox. It relies almost solely on matrix-vector multiplications,
so it can be applied to very large sparse matrices. It does not require
much memory, and it is at least as fast as a single application of
un-preconditioned LSQR to solve a least-squares problem. The method
never returns an overestimate of the condition number.

In many cases, the method is orders-of-magnitude faster than competing
methods, especially if $A$ is large and has no sparse triangular
factorization.

However, the performance of the method depends on the distribution
of the singular values of $A$, and some distributions lead to very
slow convergence. When convergence is slow or non-existent, the method
still provides a lower bound on the condition number, but it may be
loose. In such cases, methods that are based on orthogonal or triangular
factorizations or on preconditioned iterative solvers may be faster.

Our method is primarily based on one known property: that the forward
error in LSQR tends to converge to an approximate singular vector
associated with $\sigma_{\min}$. This property of LSQR and related
Krylov-subspace solvers is normally seen as a deficiency (because
it slows down the convergence to the minimizer), but it turns out
to be beneficial for condition-number estimation.

We have not explored whether a similar technique can be used in other
least-squares Krylov methods, such as LSMR~\cite{fong11}, MINRES~\cite{paige75},
and MINRES-QLP~\cite{choi11}.

\paragraph*{Acknowledgments}

Most of the work was done while Haim Avron was at IBM T.J. Watson
Research Center. We thank Nira Dyn for helpful discussions. Sivan
Toledo and Alex Druinsky were supported in part by grant 1045/09 from
the Israel Science Foundation (founded by the Israel Academy of Sciences
and Humanities) and by grant 2010231 from the US-Israel Binational
Science Foundation. Haim Avron acknowledges the support from the XDATA
program of the Defense Advanced Research Projects Agency (DARPA),
administered through Air Force Research Laboratory contract FA8750-12-C-0323.

\bibliographystyle{plain}
\bibliography{condest}

\end{document}